\documentclass[reqno]{amsart}

\usepackage{amsmath,amssymb,color}

\usepackage[latin1]{inputenc}

\newtheorem{theorem}{Theorem}[section]
\newtheorem{lemma}[theorem]{Lemma}
\newtheorem{proposition}[theorem]{Proposition}
\newtheorem{corollary}[theorem]{Corollary}
\newtheorem{remark}[theorem]{Remark}

\numberwithin{equation}{section}

\newcommand{\mc}[1]{{\mathcal #1}}
\newcommand{\mf}[1]{{\mathfrak #1}}
\newcommand{\mb}[1]{{\mathbf #1}}
\newcommand{\bb}[1]{{\mathbb #1}}

\newcommand{\<}{\langle}
\renewcommand{\>}{\rangle}

\begin{document}

\title[Exclusion processes with conductances]{Hydrodynamic limit
  of gradient exclusion processes with conductances on $\bb Z^d$}

\author[F. J. Valentim]{Fábio J. Valentim}
\address{Fábio J. Valentim \hfill\break\indent
IMPA \hfill\break\indent
Estrada Dona Castorina 110, \hfill\break\indent
J. Botânico, 22460 Rio
de Janeiro, Brasil}
\email{valentim@impa.br}

\thanks{Research supported by CNPq}

\noindent


\begin{abstract}

 Fix a smooth function $\Phi : [l,r] \to \bb R$,
  defined on some interval $[l,r]$ of $\bb R$, such that $0<b \le
  \Phi'\le b^{-1}$.
 We prove that the evolution, on the diffusive  scale, of the empirical
 density of exclusion processes in $\bb Z^d$, with  conductances given by special class of functions $W$,
 is described by the weak solutions of the
 non-linear parabolic partial differential equation
 $\partial_t \rho = \sum^d_{k=1}(d/dx_k)(d/dW_k)\Phi(\rho)$. 
  We also derive some
 properties of the operator $\sum^d_{k=1}(d/dx_k)(d/dW_k)$.
\end{abstract}

\maketitle

\section{Introduction}
\label{sec1}

We consider exclusion processes with conductances given by a special class of functions 
$W:\bb R^d\to \bb R$, 
such that
$ W(x_1,\ldots,x_d) = \sum^d_{k=1}W_k(x_k)$, 
where $d>1$ and each function $W_k: \bb R \to \bb R$ is strictly increasing, right continuous with left limits (c\`adl\`ag) , and periodic in the sense that  $W_k(u+1) - W_k(u) = W_k(1) - W_k(0)$ for all $u\in \bb R$.
 We show that, on the diffusive scale, the macroscopic evolution 
 of the empirical density of exclusion processes is described by the nonlinear 
 differential equation
\begin{equation}
\label{eq01}
\partial_t \rho \;=\; \sum^d_{k=1}\frac{d}{dx_k} \frac d{dW_k} \Phi(\rho)\;,
\end{equation}
where $\Phi$ is a smooth function strictly increasing in the range of  $\rho$ such that $0<b \le \Phi'\le b^{-1}$ and $\frac d{dW_k}$ denotes the generalized derivative, see \cite{dyn2,TC} and a revision in section \ref{lw}. In Theorem \ref{t01} we show that the operator 
$\sum^d_{k=1}(d/dx_k) (d/dW_k)$, 
defined on an appropriate domain, is non-positive, self-adjoint and dissipative; moreover, its eigenvalues are countable and have finite multiplicity, the associated eigenvectors forming a complete orthonormal system. Thus, we obtain the infinitesimal generator of a reversible Markov process.
These properties have been proved in \cite{TC} for the one-dimensional case, i.e., $d=1$.

The main tool used was the theory of energetic spaces and Friedrichs extension, 
see, for instance, \cite[chapter 5]{z}. In our case, we build the operator with the above properties by using the one-dimensional case, see section \ref{lw}.

The discrete version of the generator $\sum^d_{k=1}(d/dx_k) (d/dW_k)$ admits a  decomposition by generators of random walks with conductances also given by $W$. This allows using the method in \cite{TC,fjl} to understand the scaling limit of the process.

We consider auxiliary Markov processes associated with the empirical measure acting in the resolvent of the random walk. Since the trajectories are càdlàg, it is usual that these processes are endowed with the Skorohod topology \cite{fjl,jko,mp,kl}. However, here, the processes are endowed with the uniform topology, analogous to \cite{TC}.

As said in the introduction of \cite{TC}, non-linear versions of the partial differential equation \eqref{eq01}, $d=1$, appear naturally as scaling limits of interacting particle systems in inhomogeneous media. They may model diffusions in which permeable membranes, at the points of discontinuities of $W$, tend to reflect particles, creating space discontinuities in the solutions. But when $d\ge 2$, it is not obvious that we should also have this effect. In fact, the particles may have other options of movement. However, for this special class of functions that we are considering here, that also provides conductances for the process, we have the same effect found in the one-dimensional case.

Models with conductances have attracted the attention of several authors. An extensive list can be found at \cite{TC,s}. Recently \cite{f} has shown homogenization results for the random walk among random conductances on an infinite cluster in $\bb Z^d$. In \cite{ma}, the author proves an almost sure invariance principle for a random walker among i.i.d. conductances in $\bb Z^d$, $d\ge 2$. 

The paper is structured as follows. In Section \ref{sec2} we present the dynamics of the above exclusion process, formalize the notations used in the paper and list the main results. In Section \ref{lw} we build the operator $\sum^d_{k=1}(d/dx_k) (d/dW_k)$ with the properties listed above. Section \ref{sec04} is a preparation for Section \ref{sec3}, the discrete exclusion process is decomposed in terms of the random walk and we prove some results involving the transition and resolvent functions of the processes involved. In section \ref{sec3}, we prove the scaling limit. Finally, in Section \ref{sec4} we show that the solutions of \eqref{eq01} have finite energy.

\section{Notation and Results}
\label{sec2}

We examine the hydrodynamic behavior of a $d$-dimensional exclusion
process, $d>1$, with conductances given by a special class of functions
 $W: \bb R^d \to \bb R$ such that:
\begin{equation}
\label{w}
W(x_1,\ldots,x_d) = \sum^d_{k=1}W_k(x_k)
\end{equation}
where $W_k: \bb R \to \bb R$ are right continuous with left
limits (c\`adl\`ag) \emph{strictly increasing} functions,
 periodic in the sense that
 $$W_k(u+1) - W_k(u) = W_k(1) - W_k(0)$$ for all $u\in \bb R$ and $k=1,\ldots,d$. 
 To keep notation simple, we assume that $W_k$ vanishes at the origin,
$W_k(0)=0$.

Denote by $\bb T^d = [0, 1)^d$ the $d$-dimensional torus and by ${e_1,\ldots ,e_d}$
 the canonical basis of $\bb R^d$. For this class of functions we have:
		\begin{itemize}
		\item $W(0) = 0$,
		\item $W$ is strictly increasing on each coordinate: $$W(x + ae_j) > W(x)$$
		for all $ 1\le j\le d$, $a > 0, x \in \bb R^d$;
		\item $W$ is continuous from above: 
		$$W(x) = \lim_{y\to x,\;y\ge x}W(y),$$ 
		where we say that $y\ge x$ if $y_j\ge x_j$ for all $ 1\le j\le d$.
		\item $W$ is defined on the torus $\bb T^d$: $$W(x_1,\ldots, x_{j-1}, 0,
		x_{j+1},\ldots , x_d) = W(x_1,\ldots , x_{j-1}, 1, x_{j+1},\ldots, x_d) -
		W(e_j),$$ for all $ 1\le j\le d$, $(x_1,\ldots , x_{j-1},
		x_{j+1}, . . . , x_d) \in \bb T^{d-1}$.
		\end{itemize}

Unless explicitly stated $W$ belongs to this class.
Let $\bb T^d_N$ be the $d$-dimensional discrete torus with $N^d$
points. Distribute particles throughout $\bb T^d_N$ in such a way that each site
of $\bb T^d_N$ is occupied by at most one particle. Denote by $\eta$ the
configurations of the state space $\{0,1\}^{\bb T^d_N}$, so that $\eta(x) =0$
if site $x$ is vacant and $\eta(x)=1$ if site $x$ is occupied.

Fix $a> -1/2\;$ and $W$, for $x=(x_1,\ldots,x_d)\in \bb T^d_N$ let
\begin{equation*}
c_{x,x+e_j}(\eta) \;=\; 1 \;+\; a \{ \eta(x-e_j) + \eta(x+2\ e_j)\}\;,
\end{equation*}
where all sums are modulo $N$, and let
\begin{equation*}
\xi_{x, x+e_j} \;=\; \frac 1{N[W((x+e_j)/N) - W(x/N)]}\;=\;
\frac 1{N[W_j((x_j+1)/N) - W_j(x_j/N)]}.
\end{equation*}

The stochastic evolution can be described as follows. Let
$x=(x_1,\ldots,x_d)\in \bb T^d_N$.
At rate $\xi_{x,x+e_j}c_{x,x+e_j}(\eta)$ the occupation variables
$\eta(x)$, $\eta(x+e_j)$ are exchanged.
If $W$ is differentiable at $x/N\in[0,1]^d$, the rate at
which particles are exchanged is of order $1$ for each direction, but if some
$W_j$ is discontinuous at $x_j/N$ , the rate is of order $1/N$.
Assume, to fix ideas, that $W_j$ is discontinuous at $x_j/N$, smooth on the segment
$(x_j/N, x_j/N + \varepsilon e_j)$, $(x_j/N - \varepsilon e_j, x_j/N )$ and that $W_k$ is
differentiable in $x_k/N$ for $k \neq j$. In this case, the
rate at which particles jump over the bond $\{x-e_j, x\}$ is of order
$1/N$, while in a neighborhood of size $N$ of this bond, particles
jump at rate $1$. In particular, a particle at site $x-e_j$ jumps to $x$
at rate $1/N$ and jumps at rate $1$ to each of the $2d-1$ other options.
Particles, therefore, tend to avoid the bond $\{x-e_j,x\}$.
For the one-dimensional case (see \cite{TC}) it was shown that, on a time interval
of length $N^2$, a particle
spends a time of order $N$ at site $x$, hence particles will jump slower
over the bond $\{x-e_j, x\}$. This bond may, for instance,  model a membrane
which obstructs the passage of particles. However, in the $d$-dimensional case,
particles have the possibility  to go from $x-e_j$ to $x$, without having to jump
over the bond $\{x-e_j,x\}$. One may argue that these discontinuity points would not serve as barriers
anymore. However, for the same time interval and scaling considered in the
one-dimensional case, a particle will jump slower over the bond $\{x-e_j, x\}$. 
This is due to the fact that any
path that begins at $x-e_j$ and ends at $x$, or vice-versa, will necessarily have a $j$-th
coordinate $\{x_j -1, x_j\}$, for some $x_j$. Then, this process also models membranes
that obstruct passages of particles. Notice that these membranes are $(d-1)$-dimensional
hyperplanes embedded in a $d$-dimensional environment.

The effect of the factor $c_{x,x+e_j}(\eta)$ is analogous to the one-dimensional case.
If the parameter $a$ is positive, the presence of particles in the neighboring
sites of the bond $\{x,x+e_j\}$ speeds up the exchange rate by a factor of
order one.

The dynamics informally presented describes a Markov evolution. The
generator $L_N$ of this Markov process acts on functions $f:
\{0,1\}^{\bb T^d_N} \to \bb R$ as
\begin{equation}
\label{g4}
(L_N f) (\eta) \;=\;\sum^d_{j=1} \sum_{x \in \bb T^d_N} \xi_{x,x+e_j}\, c_{x,x+e_j}(\eta)\,
\{ f(\sigma^{x,x+e_j} \eta) - f(\eta) \} \;,
\end{equation}
where $\sigma^{x,x+e_j} \eta$ is the configuration obtained from $\eta$
by exchanging the variables $\eta(x)$ and $\eta(x+e_j)$:
\begin{equation}
\label{g5}
(\sigma^{x,x+e_j} \eta)(y) \;=\;
\begin{cases}
\eta (x+e_j) & \text{ if } y=x,\\
\eta (x) & \text{ if } y=x+e_j,\\
\eta (y) & \text{ otherwise}.
\end{cases}
\end{equation}

 A straightforward computation shows that the Bernoulli product measures
$\{\nu^N_\alpha : 0\le \alpha \le 1\}$ are invariant, and in fact
reversible, for the dynamics. The measure $\nu^N_\alpha$ is obtained
by placing a particle at each site, independently from the other
sites, with probability $\alpha$. Thus, $\nu^N_\alpha$ is a product
measure over $\{0,1\}^{\bb T^d_N}$ with marginals given by
\begin{equation*}
\nu^N_\alpha \{\eta : \eta(x) =1\} \;=\; \alpha
\end{equation*}
for $x$ in $\bb T^d_N$. For more details see \cite[chapter 2]{kl}.
We will often omit the index $N$ on $\nu^N_\alpha$.

Denote by $\{\eta_t : t\ge 0\}$ the Markov process on $\{0,1\}^{\bb
  T^d_N}$ associated to the generator $L_N$ \emph{speeded up} by
$N^2$. Let $D(\bb R_+, \{0,1\}^{\bb T^d_N})$ be the path space of
c\`adl\`ag trajectories with values in $\{0,1\}^{\bb T^d_N}$. For a
measure $\mu_N$ on $\{0,1\}^{\bb T^d_N}$, denote by $\bb P_{\mu_N}$ the
probability measure on $D(\bb R_+, \{0,1\}^{\bb T^d_N})$ induced by the
initial state $\mu_N$ and the Markov process $\{\eta_t : t\ge 0\}$.
Expectation with respect to $\bb P_{\mu_N}$ is denoted by $\bb
E_{\mu_N}$.

\subsection{The operator $\mc L_W$}

Fix $W=\sum^d_{k=1}W_k$ as in \eqref{w}.
 In \cite{TC} it is shown the existence of self-adjoint operators $\mc L_{W_k}:
 \mc D_{W_k}\subset L^2(\bb T) \to L^2(\bb T)$. Further, the set $\mc A_{W_k}$ 
 of the eigenvectors of $\mc L_{W_k}$ forms a complete orthonormal system in $L^2(\bb T)$.
  Let 
   $$\mc A_W\;=\;\{f: \bb T^d \rightarrow \bb R;f(x_1,\ldots , x_d)=
  \prod^{d}_{k=1}{ f_k(x_k)}, f_k \in  \mc A_{W_k},\;k=1,\ldots,d \}.$$

 Denote by $span(A) $ the space of finite linear combinations of the set $A$, and $\bb D_W:=span(\mc A_W)$.
  Define the operator $\bb L_W:\bb D_W \to L^2(\bb T^d)$ as follows.
  For $f =\prod^d_{k=1}{f_k}\in \mc A_W$, we have
\begin{equation}
\bb L_W(f)(x_1,\ldots x_d)=\sum^{d}_{k =1} \prod^{d}_{j=1, j\neq k} {f_j(x_j)}\mc L_{W_k}f_k(x_k),
\end{equation}
and we then extend to $\bb D_W$ by linearity. 

Lemma \ref{f17}, in Section \ref{lw}, shows that $\bb L_W$ is symmetric and non-positive;
 $\bb D_W$ is dense in $L^2(\bb T^d)$; and
the set $\mc A_W$ forms a complete, orthonormal, countable system
of eigenvectors for the operator $\bb L_W$. Let
$\mc A_W= \{ h_k\}_{ k\ge 0}$, $\{ \alpha_k\}_{ k\ge 0} $ be the corresponding eigenvalues of $-\bb L_W$, 
and consider $\mc D_W = \{v=\sum^{\infty}_{k=1}{v_kh_k}\in L^2(\bb T^d);
\sum^{\infty}_{k=1}{v^2_k\alpha^2_k}<+\infty \}$. We
define the operator $\mc L_W:\mc D_W\to L^2(\bb T^d)$ by
\begin{equation}
-\mc L_Wv = \sum^{+\infty}_{k=1}{\alpha_{k}v_{k}h_{k}}
\end{equation}

The operator $\mc L_W$ is clearly an extension of the operator $\bb L_W$,
and we present in Theorem \ref{t01} some properties of this operator.

\begin{theorem}
\label{t01}
The operator $\mc L_W : \mc D_W \to L^2(\bb T^d)$ enjoys the
following properties.

\renewcommand{\theenumi}{\alph{enumi}}
\renewcommand{\labelenumi}{{\rm (\theenumi)}}

\begin{enumerate}
\item  The domain $\mc D_W$ is
 dense in $L^2(\bb T^d)$. In particular, the set of eigenvectors
$\mc A_W=\{h_k\}_{k\ge 0}$ forms a complete orthonormal system;
\item The eigenvalues of the operator $- \mc L_W$ form a countable set
  $\{\alpha_k\}_{ k\ge 0}$. All eigenvalues have finite multiplicity,
 and it is   possible to obtain a re-enumeration $\{\alpha_k\}_{k\ge 0}$ such that
$$0= \alpha_0 \le \alpha_1 \le \cdots \;\;\;\text{and}\;\; \lim_{n\to\infty} \alpha_n  = \infty;$$

\item The operator $\bb I - \mc L_W : \mc D_W \to L^2(\bb T^d)$ is  bijective;
\item $\mc L_W: \mc D_W \to L^2(\bb T^d)$ is self-adjoint and
  non-positive:
\begin{eqnarray*}
\< -\mc L_W f , f\> \;\ge\; 0;
\end{eqnarray*}

\item $\mc L_W$ is dissipative.

\end{enumerate}
\end{theorem}

In view of (a), (b) and (d), by Hille-Yosida theorem, $\mc L_W$ is
the generator of a strongly continuous contraction
semi-group $\{P_t:L^2(\bb T^d) \to L^2(\bb
T^d) \;\}_{t\ge 0}$.

Denote by $\{G_\lambda:L^2(\bb T^d) \to L^2(\bb T^d)\;\}_{\lambda >0}$
 the semi-group of resolvents associated to the operator
$\mc L_W$: $G_\lambda = (\lambda - \mc L_W)^{-1}$. $G_\lambda$ can also be written in terms of the
semi-group $\{P_t\;;t\ge 0 \}$: $$G_\lambda = \int_0^\infty e^{-\lambda t} P_t \,
dt.$$ 

In Section \ref{sec04} we derive some properties and obtain some results
for these operators.

\subsection{The hydrodynamic equation}
\label{ss2.3}
A sequence of probability measures $\{\mu_N : N\geq 1 \}$ on $\{0,1\}^{\bb T^d_N}$
is said to be associated to a profile $\rho_0 :\bb T^d \to [0,1]$ if
\begin{equation}
\label{f09}
\lim_{N\to\infty}
\mu_N \left\{ \, \Big\vert \frac 1{N^d} \sum_{x\in\bb T^d_N} H(x/N) \eta(x)
- \int H(u) \rho_0(u) du \Big\vert > \delta \right\} \;=\; 0
\end{equation}
for every $\delta>0$ and every continuous function $H: \bb T^d \to \bb
R$. For details, see \cite[chapter 3]{kl}.

For a positive integer $m\ge 1$, denote by $C^m(\bb T^d)$ the space of
continuous functions $H:\bb T^d\to \bb R$ with $m$ continuous derivatives.
Fix $l<r$ and a smooth function $\Phi :[l,r]\to \bb R$ whose derivative 
is bounded below by a strictly positive constant and
bounded above by a finite constant:
\begin{equation*}
0 \;<\, B^{-1} \le \Phi'(x)\; \le\; B
\end{equation*}
for all $x \in [l,r]$. 
Let $\gamma : \bb T^d \to [l,r]$ be a bounded density profile and consider the parabolic differential equation
\begin{equation}
\label{g03}
\left\{
\begin{array}{l}
{\displaystyle \partial_t \rho \; =\; \mc L_W \Phi(\rho) } \\
{\displaystyle \rho(0,\cdot) \;=\; \gamma(\cdot)}
\end{array}
\right. .
\end{equation}

A bounded function $\rho : \bb R_+ \times \bb T^d \to [l,r]$
is said to be a weak solution of the parabolic differential equation \eqref{g03} if
\begin{equation*}
\< \rho_t, G_\lambda H\> \;-\; \< \gamma , G_\lambda H\>
\;=\; \int_0^t \< \Phi(\rho_s) , \mc L_W G_\lambda H \>\, ds\;
\end{equation*}
for every continuous function $H:\bb T^d\to \bb R$, all $t>0$ and all $\lambda >0$.

Existence follows from tightness of the sequence of probability measures
$\bb Q_{\mu_N}^{W,N}$ introduced in Section \ref{sec3}.
The proof of uniqueness of weak solutions is analogous to \cite{TC}.

\begin{theorem}
\label{t02}
Fix a continuous initial profile $\rho_0 : \bb T^d \to [0,1]$ and
consider a sequence of probability measures $\mu_N$ on $\{0,1\}^{\bb
  T^d_N}$ associated to $\rho_0$, in the sense of \ref{f09}. Then, for any $t\ge 0$,
\begin{equation*}
\lim_{N\to\infty}
\bb P_{\mu_N} \left\{ \, \Big\vert \frac 1{N^d} \sum_{x\in\bb T^d_N}
H(x/N) \eta_t(x) - \int H(u) \rho(t,u)\, du \Big\vert
> \delta \right\} \;=\; 0
\end{equation*}
for every $\delta>0$ and every continuous function $H$. Here, $\rho$
is the unique weak solution of the non-linear equation \eqref{g03}
with $l=0$, $r=1$, $\gamma = \rho_0$ and $\Phi(\alpha) = \alpha + a
\alpha^2$.
\end{theorem}

\begin{remark}
\label{sf01}
The specific form of the rates $c_{x,x+e_i}$ is not important, but two
conditions must be fulfilled. The rates must be strictly positive,
they may not depend on the occupation variables $\eta(x)$, $\eta(x+e_i)$,
but they have to be chosen in such a way that the resulting process is
\emph{gradient}. (cf. Chapter 7 in \cite{kl} for the definition of
gradient processes).

We may define rates $c_{x,x+e_i}$ to obtain any polynomial $\Phi$ of the
form $\Phi(\alpha) = \alpha + \sum_{2\le j\le m} a_j \alpha^j$, $m\ge
1$, with $1+ \sum_{2\le j\le m} j a_j >0$. Let, for instance, $m=3$
then the rates
\begin{align*}
\hat c_{x,x+e_i} (\eta)\;\;& =\;\; c_{x,x+e_i} (\eta)\;\;+\\
&b\left\{ \eta(x-2e_i) \eta(x-e_i) + \eta(x-e_i) \eta(x+2e_i) + \eta(x+2e_i)
\eta(x+3e_i)\right\},
\end{align*}
satisfy the above three conditions, 
where $c_{x,x+e_i}$ is the rate defined at the beginning of Section 2 and
$a$, $b$ are such that $1+2a + 3b>0$. An elementary computation shows
that  
$\Phi(\alpha) = 1 + a \alpha^2 + b \alpha^3$.
\end{remark}

In Section \ref{sec4} we prove that any limit point $\bb Q^*_{W}$ of the
sequence $\bb Q_{\mu_N}^{W,N}$ is concentrated on trajectories
$\rho(t,u) du$ with finite energy in the following sense:
for each $1\le j\le d$, there is a Hilbert space $L^2_{x_j\otimes W_j}$,
associated to $W_j$, such that
$$\int^t_0 ds\;\|\frac d {dW_j}\Phi(\rho(s,.))\|_{x_j\otimes W_j}^2<\infty\;,$$
where $\|.\|_{x_j\otimes W_j}$ is the norm in $L^2_{x_j\otimes W_j}$ and $d/dW_j$
is the derivative, which must be understood in the generalized sense.

\begin{section}{The operator $\mc L_W$}
\label{lw}

The operator $\mc L_W:\mc D_W\subset L^2(\bb T^d) \rightarrow L^2(\bb T^d)$
is a natural extension, for the $d$-dimensional case, of the self-adjoint operator obtained
for the one-dimensional case in \cite{TC}.

We begin by presenting one of the main results obtained in \cite{TC}, and we then 
present the necessary modifications to conclude similar results for the $d$-dimensional case.

\subsection{Some remarks on the one-dimensional case }
Denote by $\<\cdot , \cdot\>$ the inner product of $L^2(\bb T),$
where $\bb T\subset \bb R$ is the one-dimensional torus:
\begin{equation*}
\< f, g\>\;=\; \int_{\bb T} f(u)\, g(u)\, du\;.
\end{equation*}

Let $W_1: \bb R \to \bb R$ be a \emph{strictly increasing} right continuous
function, with left limits (c\`adl\`ag), and periodic in the sense that
 $W_1(u+1) - W_1(u) = W_1(1) - W_1(0)$ for all $u$ in $\bb R$.

Let $\mc D_{W_1}$ be the set of functions $f$ in $L^2(\bb T)$ such that
\begin{equation*}
f(x) \;=\; a \;+\; b W_1(x) \;+\; \int_{(0,x]} W_1(dy) \int_0^y \mf f(z) \, dz,
\end{equation*}
for some function $\mf f$ in $L^2(\bb T)$ such that
\begin{equation*}
\int_0^1 \mf f(z) \, dz \;=\; 0\;, \quad
\int_{(0,1]} W_1(dy) \Big ( b + \int_0^y \mf f(z) \, dz \Big) \;=\;0\;
\end{equation*}
Define the operator $\mc L_{W_1} : \mc D_{W_1} \to L^2(\bb T)$ by $\mc
L_{W_1} f = \mf f$. Formally
\begin{equation}
\label{c01}
\mc L_{W_1} f \;=\; \frac d{dx} \, \frac d{dW_1} \, f\;,
\end{equation}
where the generalized derivative $d/dW_1$ is defined as
\begin{equation}
\label{c02}
\frac{d f}{dW_1} (x) = \lim_{\epsilon\rightarrow 0} \frac{f(x+\epsilon)
-f(x)}{W_1(x+\epsilon) -W_1(x)}\;,
\end{equation}
if the above limit exists and is finite.

Denote by $\bb I$ the identity operator in $L^2(\bb T)$.

\begin{theorem}
\label{s11}
The operator $\mc L_{W_1} : \mc D_{W_1} \to L^2(\bb T)$ enjoys the
following properties:

\renewcommand{\theenumi}{\alph{enumi}}
\renewcommand{\labelenumi}{{\rm (\theenumi)}}

\begin{enumerate}
\item $\mc D_{W_1}$ is dense in $L^2(\bb T)$;

\item The operator $\bb I - \mc L_{W_1} : \mc D_{W_1} \to L^2(\bb T)$ is
  bijective;

\item $\mc L_{W_1}: \mc D_{W_1} \to L^2(\bb T)$ is self-adjoint and
  non-positive:
\begin{eqnarray*}
\< -\mc L_{W_1} f , f\> \;\ge\; 0;
\end{eqnarray*}

\item $\mc L_{W_1}$ is dissipative i.e., for all $g\in \mc D_W$ and $\lambda>0$, we have
\begin{equation*}
\|\lambda g \|\le \|(\lambda \bb I-\mc L_{W_1})g\|;
\end{equation*}

\item The eigenvalues of the operator $- \mc L_W$ form a countable set
  $\{\lambda_n : n\ge 0\}$. All eigenvalues have finite multiplicity,
  $0= \lambda_0 \le \lambda_1 \le \cdots$, and $\lim_{n\to\infty} \lambda_n
  = \infty$;

\item The eigenvectors $\{f_n\}_{n\ge0}$ of the operator $\mc L_W$ form a complete orthonormal system.

\end{enumerate}
\end{theorem}

The proof can be found in \cite{TC}.

\subsection{The $d$-dimensional case}

 Consider $W$ as in \eqref{w}. Let $\mc A_{W_k}$ be the countable complete orthonormal
 system of eigenvectors for the operator
  $\mc L_{W_k}: \mc D_{W_k}\subset L^2(\bb T) \rightarrow \bb R$ given by
Theorem \ref{s11}. Let 
   $$\mc A_W\;=\;\{f: \bb T^d \rightarrow \bb R;f(x_1,\ldots , x_d)=
  \prod^{d}_{k=1}{ f_k(x_k)}, f_k \in  \mc A_{W_k } \}.$$

Note that, by  Fubini's  theorem, the set $\mc A_W$ is orthonormal in
$L^2(\bb T^d)$, and the constant functions are eigenvectors for the operator
$\mc L_{W_k}$. Moreover, $\mc A_{W_k}\subset \mc A_W$, in the sense that
$f_k(x_1,\ldots ,x_d) =f_k(x_k),\;\;f_k\in \mc A_{W_k}$.

Define the operator $\bb L_W:\bb D_W:=span(\mc A_W) \to L^2(\bb T^d)$ as follows: 
for $f =\prod^d_{k=1}{f_k}\in \mc A_W$,
\begin{equation}
\label{eq31}
\bb L_W(f)(x_1,\ldots x_d)=\sum^{d}_{k =1} \prod^{d}_{j=1, j\neq k}
{f_j(x_j)}\mc L_{W_k}f_k(x_k),
\end{equation}
and extend to $\bb D_W$ by linearity.

By \eqref{c01}, the operators $\mc L_{W_k}$ can be formally 	
extended to functions defined on $\bb T^d$ as follows.
Given a function $f:\bb T^d\to \bb R$, we define $\mc L_{W_k}f$ as 
\begin{equation}
\label{f008}
\mc L_{W_k}f\;=\; \frac d{dx_k} \, \frac d{dW_k} \, f,\;
\end{equation}
where the generalized derivative $d/dW_k$ is defined by
\begin{equation}
\label{f004}
\frac{d f}{dW_k} (x_1,\ldots,x_k,\ldots, x_d) = \lim_{\epsilon\rightarrow 0}
\frac{f(x_1,\ldots,x_k +\epsilon,\ldots, x_d)
-f(x_1,\ldots,x_k,\ldots, x_d)}{W_k(x_k+\epsilon) -W_k(x_k)}\;,
\end{equation}
if the above limit exists and is finite. Hence, by \eqref{eq31}, if $f\in \bb D_W$ 

\begin{equation}
\label{f002}
\bb L_W f=\sum^{d}_{k =1}\mc L_{W_k}f.
\end{equation}

Note that if $f=\prod^d_{k=1}{f_k}$, where $f_k\in \mc A_{W_k}$ is an eigenvector
of $\mc L_{W_k}$ associated to the eigenvalue $\lambda _k$, then $f$ is an eigenvector of
 $\bb L_W$ with eigenvalue $\sum^d_{k=1}{\lambda _k}$. 

\begin{lemma}
\label{f17}
The following statements hold:

\renewcommand{\theenumi}{\alph{enumi}}
\renewcommand{\labelenumi}{{\rm (\theenumi)}}

\begin{enumerate}
\item The set $\bb D_W$ is dense in $L^2(\bb T^d)$;
\item The operator $\bb L_W : \bb D_W \to L^2(\bb T^d)$ is symmetric and
  non-positive:
\begin{eqnarray*}
\< -\bb L_W f , f\> \;\ge\; 0.
\end{eqnarray*}
\end{enumerate}

\end{lemma}

\begin{proof}
The strategy to prove the above lemma is the following. We begin by showing that the set
$$\mc S\; =\; span(\{f\in L^2(\bb T^d) ; f(x_1,\ldots , x_d) = \prod^d_{k=1}{f_k(x_k)}, f_k\in \mc D_{W_k}\})$$ 
is dense in 
$$\bb S = span(\{f\in L^2(\bb T^d) ; f(x_1,\ldots , x_d) =  \prod^d_{k=1}f_k(x_k),f_k\in L^2(\bb T)\}).$$

We then show that $\bb D_W$ is dense in $\mc S$.
Since $\bb S$ is dense in $L^2(\bb T^d)$, item (a) follows.

 We now prove item (a) rigorously. Since $\mc S$ is a vector space, we only have to show that we can approximate the functions
$\prod^d_{k=1}f_k\in L^2(\bb T^d)$, where $f_k\in \mc D_{W_k}$, by functions of $\bb D_W$. 
By Theorem \ref{s11}, the set $\mc D_{W_k}$
 is dense in  $L^2(\bb T)$, thus, there is a sequence $(f^k_n)_{n\in \bb N}$
 converging to $f_k$ in $L^2(\bb T)$. Thus
 Let  $$f_n(x_1,\ldots , x_d) = \prod^d_{k=1}{f^k_n(x_k)}.$$
 By the triangle inequality and Fubini's theorem, the sequence $(f_n)$ converges to
 $\prod^d_{k=1}f_k$. Fix $\epsilon > 0$ and let
 $$h(x_1,\ldots ,x_d)=\prod^d_{k=1}{h_k(x_k)}, \;\;\;\;h_k 
 \in \mc D_{W_k}. $$

 Since for each $k=1\ldots,d,\;\; \mc A_{W_k}\subset \mc D_{W_k}$ is a complete
 orthonormal set, there exist sequences $g^k_j \in \mc A_{W_k}$ and
 $\alpha ^k_j\in \bb R$ such that
 $$\|h_k - \sum^{n(k)}_{j=1}{\alpha ^k_jg^k_j}\|_{L^2(\bb T)} < \delta \;,$$
 where $\delta = \epsilon/{dM^{d-1}}$ and $M:=1+ sup_{k=1:n}\left\| h_k\right\| $.
Let $$g(x_1,\ldots ,x_d) = \prod^d_{k=1}{\sum^{n(k)}_{j=1}{\alpha^k_jg^k_j(x_k)}}\in \bb D_W.$$

An application of the triangle inequality and Fubini's theorem yields
$\left\|h-g\right\|<\epsilon$. This proves (a).
 
 To prove (b), let
\begin{eqnarray*}
f(x_1,\ldots , x_d)=\prod^d_{k=1}{f_k(x_k)}\;\;\text{and}\;\;
g(x_1,\ldots , x_d)=\prod^d_{k=1}{g_k(x_k)}
\end{eqnarray*}
be functions belonging to $\mc A_W$. We have that
\begin{equation*}
\<f,\bb L_Wg\> = \<\prod^d_{k=1}{f_k},\sum^d_{k=1}
{\prod^d_{j=1,j\neq k}{g_j}\mc L_{W_k}g_k}\> =
\sum^d_{k=1}\<\prod^d_{j=1,j\neq k}{f_jg_j},f_k\mc L_{W_k}g_k\>,
\end{equation*}
where $\<\cdot,\cdot\>$ denotes the inner product in $L^2(\bb T^d)$.
Since, by Theorem \ref{s11}, $\mc L_{W_k}$ is self-adjoint, we have
\begin{equation*}
\sum^d_{k=1}\<\prod^d_{j=1,j\neq k}{f_jg_j},g_k\mc L_{W_k}f_k\> = \<\mc L_Wf,g\>.
\end{equation*}

In particular, the operator $\mc L_{W_k}$
is non-positive and therefore
\begin{equation*}
\<f,\bb L_Wf\> = \sum^d_{k=1}\<\prod^d_{j=1,j\neq k}{f^2_j},f_k\mc L_{W_k}f_k\>\le0.
\end{equation*}
Item (b) follows by linearity.
\end{proof}

Lemma \ref{f17} implies that the set $\mc A_W$ forms a complete orthonormal
 countable system of eigenvectors for the operator
$\bb L_W$. Let $\mc A_W= \{ h_k\}_{ k\ge 0}$, and let $\{ \alpha_k\}_{ k\ge 0} $ be
the corresponding eigenvalues of $-\bb L_W$. Consider
$$\mc D_W = \{v=\sum^{\infty}_{k=1}{v_kh_k}\in L^2(\bb T^d);
\sum^{\infty}_{k=1}{v^2_k\alpha^2_k}<+\infty \},$$
and define

\begin{equation}
-\mc L_{W}v = \sum^{+\infty}_{k=1}{\alpha_{k}v_{k}h_{k}}.
\end{equation}

The operator $\mc L_W$ is clearly an extension of the operator $\bb L_W$.
Formally, by \eqref{f002},
\begin{equation}
\label{f07}
\mc L_W f \; =\; \sum^d_{k=1}\mc L_{W_k}f,
\end{equation}
where
\begin{equation*}
\mc L_{W_k}f\;=\; \frac d{dx_k} \, \frac d{dW_k} \, f.
\end{equation*}

We are now in conditions to prove Theorem \ref{t01}.
\begin{proof}[Proof of Theorem \ref{t01}]
Since $\bb D_W \subset \mc D_W$, the density of $\mc D_W$ in $L^2(\bb T^d)$
follows from the density of $\bb D_W$, shown in lemma \ref{f17}.

If $\alpha_k$ are eigenvalues of $-\mc L_W$, we may find 
eigenvalues $\;\lambda_j$, associated to some $f_j \in \mc A_{W_j}$, 
such that $\alpha_k = \sum^d_{j=1}{\lambda_j}$.
By Theorem \ref{s11}[item (e)], (b) follows.

Let $\{\alpha_k\}_{k\ge0} $ be the set of eigenvalues of $\bb -\mc L_W$. 
Then the  set of eigenvalues of $\bb I-\mc L_W$ is $\{\gamma_k\}_{k\ge0} $, where
$\gamma_k=\alpha_k+1$, and the eigenvectors are the same as the ones of $\mc L_W$.
By item (b), we have
$$1= \gamma_0 \le \gamma_1 \le \cdots\;\; \text{and}\;\;\lim_{n\to\infty}
\gamma_n  = \infty \;.$$

Thus, $\bb I-\mc L_W$ is injective, and for
$$v=\sum^{+\infty}_{k=1}{v_kh_k}\in L^2(\bb T^d)\;,\;\;
\text{such that}\;\;\sum^{\infty}_{k=1}{v^2_k}<+\infty \;, $$
let
$$u= \sum^{+\infty}_{k=1}{\frac{v_k}{\gamma_k}h_k}\;,$$
then $u\in \mc D_W$ and $(\bb I-\mc L_W)u=v$. Hence, item (c) follows.

Let $\mc L^*_W:\mc D_{W^*}\subset L^2(\bb T^d) \rightarrow L^2(\bb T^d)$
be the adjoint of $\mc L_W$. Since $\mc L_W$ is symmetric, we have
$\mc D_W \subset \mc D_{W^*}$. So,  to show the equality of the operators
it suffices to show that $\mc D_{W^*} \subset \mc D_W $. Given
$$\varphi=\sum^{+\infty}_{k=1}{\varphi_k h_k} \in \mc D_{W^*,}$$ let
$\mc L_{W*}\varphi = \psi \in L^2(\bb T^d)$. Therefore, for all
$v=\sum^{+\infty}_{k=1}{v_kh_k}\in \mc D_W$,
\begin{equation*}
\< v , \psi\> = \< v ,\mc L_{W*}  \varphi \> =  \< \mc L_W v ,
\varphi\> = \sum^{+\infty}_{k=1}{-\alpha_kv_k\varphi_k}.
\end{equation*}

Hence $$\psi = \sum^{+\infty}_{k=1}{-\alpha_k\varphi_kh_k}\;,$$ in particular,
$$\sum^{+\infty}_{k=1}\alpha^2_k\varphi^2_k<+\infty \;\; \text{and}\;\;
\varphi\in \mc D_W.$$

Thus, $\mc L_W$ is self-adjoint. Let $v=\sum^{+\infty}_{k=1}{v_kh_k}\in \mc D_W$. 
 From item (b) $\alpha_k\ge0$ and
$$\<-\mc L_Wv,v\> = \sum^{+\infty}_{k=1}{\alpha_kv^2_k} \ge0.$$

Therefore $\mc L_W$ is non-positive and item (d) follows.
 
Fix a function $g$ in
$\mc D_W$, $\lambda >0$, and let $f= (\lambda \bb I - \mc L_W) g$.
Taking inner product, with respect
to $g$, on both sides of this equation, we obtain
\begin{eqnarray*}
\lambda \< g , g\> \;+\; \< -\mc L_W g , g\>
\;=\; \< g , f\> \;\le\;  \< g , g\>^{1/2} \,
\<f , f\>^{1/2}\;.
\end{eqnarray*}
Since $g$ belongs to $\mc D_W$, by (d), the second term on the left
hand side is non-negative. Thus, $\Vert \lambda g\Vert \le \Vert f \Vert =
\Vert (\lambda \bb I - \mc L_W) g \Vert $.
\end{proof}
\end{section}

\begin{section}{Random walk with conductances}
\label{sec04}
Recall the decomposition obtained in \eqref{f07}: $$\mc L_W = \sum^d_{j=1}\mc L_{W_j}.$$

The discrete version of $\mc L_W$ is the generator $L_N$ of a Markov process 
given as follows. For each function $f:\{0,1\}^{\bb T^d_N}\to \bb R$,
\begin{equation}
\label{ff}
L_Nf(\eta) \; = \; \sum^d_{j=1}L^j_Nf(\eta),
\end{equation}
where $$L^j_Nf(\eta)\;=\sum_{x\in \bb T^d_N}\xi_{x,x+e_j}c_{x,x+e_j}
(\eta)[f(\eta^{x,x+e_j})-f(\eta)].$$

In Section \ref{sec2}, we described, informally, the dynamics of this Markov evolution.

\subsection{Discrete approximation of the  operator $\mc L_W$}

Let $j=1,\ldots ,d$. Consider the random walks $\{X_t^j\;\}_{t\ge0}$ 
on the discrete torus, $N^{-1}\bb T_N$, which jumps from
$x/N$ (resp. $(x+1)/N$) to $(x+1)/N$ (resp. $x/N$) with rate
$$N^2 \xi^j_{x,x+1} =N/\{W_j((x+1)/N) - W_j(x/N)\}.$$

Let $\{X_t=(X^1_t, \ldots , X^d_t)\}_{t\ge 0}$ be the random walk
on $N^{-1}\bb T^d_N\;,$ where $\bb T^d_N$ is the discrete
$d$-dimensional torus with $N^d$ points.

The generator $\bb L_N$ of this Markov process acts on functions $f:\bb T^d_N\to \bb R$ as
\begin{equation}
\label{fg}
\bb L_N f(x/N) \; =\;\sum^d_{j=1}\bb L^j_N f(x/N),
\end{equation}
where
\begin{eqnarray*}
\bb L^j_N f(x/N) &=& N^2\big \{ \xi_{x,x+e_j} [f((x+e_j)/N) - f(x/N)] \\
&+& \xi_{x-e_j,x} [f((x-e_j)/N) - f(x/N)] \big \}
\end{eqnarray*}
are the generators of the one-dimensional random walks $\{X^j_t\}_{t\ge0}$. 

Note that $\bb L^j_N f(x/N)$ is in fact a discrete version 
of the operator $\mc L_{W_j}$. The counting measure $m_N$ on $N^{-1} \bb T^d_N$ is
reversible for this process.
\subsection{Semigroups and resolvents.}
In this subsection we introduce families of semigroups and resolvents,
associated to the generators $\bb L_N$  and $\mc L_W$.
We present some properties and results regarding the convergence of these operators.

  Denote by $\{P^N_t : t\ge 0\}$ (resp. $\{G_\lambda^N :
\lambda >0\}$) the semigroup (resp. the resolvent) associated to the
generator $\bb L_N$, by $\{P^{N,j}_t : t\ge 0\}$ the semigroup associated to the
generator $\bb L^j_N$, by $\{P^j_t :t\ge 0\}$ the semigroup associated to the generator $\mc L_{W_j}$
 and by $\{P_t : t\ge 0\}$ (resp. $\{G_\lambda:
\lambda >0\}$) the semigroup (resp. the resolvent) associated to the
generator $\mc L_W$.

Since the jump rates from
$x/N$ (resp. $(x+e_j)/N$) to $(x+e_j)/N$ (resp. $x/N$) are equal, 
$P^N_t$ is symmetric: $P^N_t(x,y)\;=\;P^N_t(y,x)$.

Using the decompositions \eqref{fg} and \eqref{f07}, we have that

\begin{equation*}
P^N_t(x,y) = \prod^d_{j=1}P^{N,j}_t(x_j,y_j)\;\;\text{and}\;\;P_t(x,y) =
\prod^d_{j=1}P^{j}_t(x_j,y_j).
\end{equation*}

By definition, for every $H: N^{-1}\bb T^d_N \to \bb R$,
\begin{equation*}
G_\lambda H  \;=\; \int_0^\infty dt\, e^{-\lambda t}  P_t H \; = \;
(\lambda\bb I-\mc L_W)^{-1}H,
\end{equation*}
where $\bb I$ is the identity operator.

\begin{lemma}
\label{l1}
Let $H:\bb T^d \rightarrow \bb R$ be a continuous function. Then
\begin{equation}
\lim_{N\to +\infty} {\frac{1}{N^d}\sum_{x\in \bb T^d_N}|P^N_tH(x/N)-P_tH(x/N)|}=0.
\end{equation}
\end{lemma}
\begin{proof}
If $H:\bb T^d \rightarrow \bb R$ has the form $H(x_1,\ldots , x_d) = \prod^d_{j=1}H_j(x_j)$, we have
\begin{equation}
\label{pt}
P^N_tH(x) \;=\; \prod^d_{j=1}P^{N,j}_tH_j(x_j)\;\;\text{and}\;\; P_tH(x) \;=\;
\prod^d_{j=1}P^{j}_tH_j(x_j).
\end{equation}

Now, for any continuous function $H:\bb T^d \rightarrow \bb R$, and any $\epsilon>0$, 
we can find continuous functions $H_{j,k}:\bb T\to\bb R$, such that $H':\bb T^d \rightarrow \bb R$ given by
 $$H'(x)= \sum^m_{j=1}\prod^d_{k=1}H_{j,k}(x_k)$$ satisfies
$\|H'-H\|_{\infty}\le\epsilon.$ Thus,
\begin{equation*}
\frac{1}{N^d}\sum_{x\in \bb T^d_N}|P^N_tH(x/N)-P_tH(x/N)|\;\le\; 2\epsilon +
\frac{1}{N^d}\sum_{x\in \bb T^d_N}|P^N_tH'(x/N)-P_tH'(x/N)|.
\end{equation*}

By \eqref{pt} and similar identities for $P_tH'$ and $P_t^{N,j}H'$, 
the sum on the right hand side in the previous inequality is less than or equal to
\begin{eqnarray*}
\frac{1}{N^d}\sum_{x\in \bb T^d_N}\sum^{m}_{j=1}|\prod^d_{k=1}P^{N,k}_tH_{j,k}(x_k/N)-
\prod^d_{k=1}P^k_tH_{j,k}(x_k/N)|\;\le\;
\\
\frac{1}{N^d}\sum_{x\in \bb T^d_N}\sum^{m}_{j=1}C_j
\sum^d_{k=1}|P^{N,k}_tH_{j,k}(x_k/N)-P^k_tH_{j,k}(x_k/N)|,
\end{eqnarray*}
where $C_j$ is a constant that depends on the product
$\prod^d_{k=1}H_{j,k}$. The previous expressions can be rewritten as
\begin{eqnarray*}
\sum^{m}_{j=1}C_j\sum^d_{k=1}\frac{1}{N^d}
\sum_{x\in \bb T^{d-1}_N}\sum^N_{i=1}|P^{N,k}_tH_{j,k}(i/N)-P^k_tH_{j,k}(i/N)|\;\;=\;
\\
\sum^{m}_{j=1}C_j\sum^d_{k=1}\frac{1}{N}
\sum^N_{i=1}|P^{N,k}_tH_{j,k}(i/N)-P^k_tH_{j,k}(i/N)|.
\end{eqnarray*}
Moreover, by \cite[lemma 4.5 item iii]{fjl}, when $N\to \infty$, the last expression converges to $0$.

\end{proof}
\begin{corollary}
\label{c1}
Let $H:\bb T^d\to \bb R$ be a continuous function. Then
\begin{equation}
\lim_{N\to +\infty}{\frac{1}{N^d}\sum_{x\in \bb T^d_N}|G^N_{\lambda}H(x/N)-G_{\lambda}H(x/N)}|=0.
\end{equation}
\end{corollary}
\begin{proof}
By definition of resolvent, for each $N$, the previous expression is less than or equal to
\begin{equation*}
\int^{\infty}_0{dt\ e^{-\lambda t}\frac 1 {N^d}\sum_{x\in\bb T^d_N}|P^N_tH(x/N)-P_tH(x/N)|}.
\end{equation*}
Corollary now follows from the previous Lemma.
\end{proof}

Let $f_N:\bb T^d_N\to\bb R$ be any function. Then, whenever needed, 
we consider $f:\bb T^d\to \bb R$ an extension of $f_N$ 
to $\bb T^d$ given by:
\begin{equation*}
f(y) = f_N(x),\; \;\text{if} \;\;x\in\bb T^d_N,\;\;y\ge x \;\;\text{and} \;\; \|y-x\|_{\infty}< \frac1N.
\end{equation*}

\begin{lemma}
\label{l03}
Let $H:\bb T^d \rightarrow \bb R$ be a continuous function.
then the extension of $P^N_tH:\bb T^d_N\to \bb R$ to $\bb T^d$ 
belongs to $L^1(\bb T^d)$, and  
$$\int_{\bb T^d}duP^N_tH(u)\;=\;\frac 1 {N^d}\sum_{x\in \bb T^d} H(x/N).$$
\end{lemma}
\begin{proof}
Assume, without loss of generality, that $H\ge 0$. Since the transition probability
$P^N_t(x, y)$ is symmetric, we have
\begin{eqnarray*}
\int_{\bb T^d}{du P^N_tH(u)}\;\;=\;\;\frac 1 {N^d}\sum_{x,y\in \bb T^d_N}P^N_t(x, y)H(y/N)\;=\;\\
\frac 1 {N^d}\sum_{y\in \bb T^d_N}H(y/N)\sum_{x\in \bb T^d_N}P^N_t(y,x)\;=\;\frac 1 {N^d}\sum_{y\in \bb T^d_N}H(y/N).
\end{eqnarray*}
This proves the identity and also that $P^N_tH\in L^1(\bb T^d)$.
\end{proof}

The next lemma shows that $H$ can be approximated by $P^N_tH$. As an
immediate consequence we obtain an approximation result involving the resolvent.
\begin{lemma}
\label{2.8}
Let $H:\bb T^d \rightarrow \bb R$ be a continuous function, then,
\begin{equation}
\lim_{t\to 0}\limsup_{N\to +\infty} {\frac{1}{N^d}\sum_{x\in \bb T^d_N}|P^N_tH(x/N)-H(x/N)|}=0,
\end{equation}
and
\begin{equation}
\lim_{\lambda\to +\infty}\limsup_{N\to +\infty} {\frac{1}{N^d}\sum_{x\in \bb T^d_N}|\lambda R^N_{\lambda}H(x/N)-H(x/N)|}=0.
\end{equation}
\end{lemma}
\begin{proof}
Fix $\epsilon> 0$ and consider $H'$ as in the proof of Lemma \ref{l1}. Thus,
\begin{equation*}
\frac{1}{N^d}\sum_{x\in \bb T^d_N}|P^N_tH(x/N)-H(x/N)|\;\le\; 2\epsilon +
\frac{1}{N^d}\sum_{x\in \bb T^d_N}|P^N_tH'(x/N)-H'(x/N)|,
\end{equation*}
where the second term on the right hand side is less than or equal to
\begin{equation*}
C_0 \sup_{j,k}\frac{1}{N^d}\sum_{x\in \bb T^d_N}|P^{N,k}_tH_{j,k}(x_k/N)-H_{j,k}(x_k/N)|,
\end{equation*}
$C_0$ being a constant that depends on $H'$.
By \cite[lemma 4.6]{fjl}, the last expression converges to $0$, when $N\to \infty$, and  then $t\to 0$.
This proves the first equality.

To obtain the second limit, note that, by definition of the resolvent,
the second expression is less than or equal to
$$ \int^{\infty}_0{dt\lambda e^{-\lambda t}\frac 1 {N^d}\sum_{x\in \bb T^d_N }|P^N_tH(x/N)- H(x/N)|}.$$

By lemma \ref{l03}, the sum is uniformly bounded in $t$ and $N$. 
By the first part of Lemma \ref{l03}, it vanishes as $N \to \infty$ and $t\to 0$.
This proves the second part.

\end{proof}

Fix a function $H:\bb T^d_N\to \bb R$. For $\lambda>0$, let $H^N_\lambda = G^N_\lambda H$ 
be the solution of the resolvent equation
$$\lambda H^N_{\lambda}-\bb L_NH^N_{\lambda}\;=\;H.$$ 

Taking inner product on both sides of this equation with respect to $H^N_{\lambda}$,we obtain
$$
\lambda \frac{1}{N^d}\sum_{x\in \bb T^d_N}(H^N_{\lambda}(x/N))^2\;-\;
\frac{1}{N^d}\sum_{x\in \bb T^d_N}H^N_{\lambda}(x/N)\bb L_NH^N_{\lambda}
$$
$$=\frac{1}{N^d}\sum_{x\in \bb T^d_N}H^N_{\lambda}(x/N)H(x/N).$$

A simple computation shows that the second term on the left hand side is equal to
\begin{equation*}
\frac{1}{N^d}\sum^d_{j=1}\sum_{x\in \bb T^d_N}\xi_{x,x+e_j}[\nabla_{N,j}H^N_{\lambda}(x/N)]^2,
\end{equation*}
where $\nabla_{N,j}H(x/N)= N[H((x+e_j)/N) - H(x/N)]$ is the discrete derivative of the function
$H$ in the direction of the vector $e_j$. In particular, by Schwarz inequality,
\begin{equation}
\label{f05}
\begin{array}{l}
{\displaystyle
\frac{1}{N^d}\sum_{x \in \bb T^d_N}
H_\lambda^N (x/N)^2 \;\leq\; \frac{1}{\lambda^2}
\frac{1}{N^d}\sum_{x \in \bb T^d_N}  H (x/N)^2 \;\;\text{and}\;\;}\\
{\displaystyle
\quad \frac{1}{N^d}\sum^d_{j=1}\sum_{x \in \bb T^d_N} \xi_{x,x+e_j}[\nabla_{N,j}
H_\lambda^N(x/N)]^2 \;\leq\; \frac{1}{\lambda}
\frac{1}{N^d}\sum_{x \in \bb T_N^d}  H (x/N)^2 \;.}
\end{array}
\end{equation}

We have proved the following.

\begin{proposition}
Let $H:\bb T^d \to \bb R$ be a continuous function, $H^N_\lambda = G^N_\lambda H$,
 and let $m_N$ be the counting measure on $\bb T^d_N$. Then
$H_\lambda^N$ and $\nabla_{N,j}H_\lambda^N$ converge to $0$ when 
$N \to\infty$ and $\lambda \to \infty$ in $L^2(\bb T^d_N,m_N)$.
\end{proposition}
\end{section}

\section{Scaling limit}
\label{sec3}

Let $\mc M$ be the space of positive measures on $\bb T^d$ with total
mass bounded by one endowed with the weak topology. Recall that
$\pi^{N}_{t} \in \mc M$ stands for the empirical measure at time $t$.
This is the measure on $\bb T^d$ obtained by rescaling space by $N$ and
by assigning mass $1/N^d$ to each particle:
\begin{equation}
\label{f01}
\pi^{N}_{t} \;=\; \frac{1}{N^d} \sum _{x\in \bb T^d_N} \eta_t (x)\,
\delta_{x/N}\;,
\end{equation}
where $\delta_u$ is the Dirac measure concentrated on $u$. 

For a continuous function $H:\bb T^d \to \bb R$, $\<\pi^N_t, H\>$ stands for
the integral of $H$ with respect to $\pi^N_t$:
\begin{equation*}
\<\pi^N_t, H\> \;=\; \frac 1{N^d} \sum_{x\in\bb T^d_N}
H (x/N) \eta_t(x)\;.
\end{equation*}
This notation is not to be mistaken with the inner product in
$L^2(\bb T^d)$ introduced earlier. Also, when $\pi_t$ has a density
$\rho$, $\pi(t,du) = \rho(t,u) du$, we sometimes write $\<\rho_t, H\>$
for $\<\pi_t, H\>$.

For a local function $g: \{0,1\}^{\bb Z^d} \to \bb R$, let
$\tilde g :[0,1]\to \bb R$ be the expected value of $g$ under the
stationary states:
\begin{equation*}
\tilde g (\alpha) \;=\; E_{\nu_\alpha} [ g(\eta)]\;.
\end{equation*}

For $\ell \ge 1$ and $d$-dimensional integer $x=(x_1,\ldots,x_d)$, denote 
by $\eta^{\ell} (x)$ the empirical density of particles in the box 
$ \bb B_+^\ell(x)= \{(y_1,\ldots,y_d)\in\bb Z^d\;;0\le y_i-x_i < \ell\}$: 
\begin{equation*}
  \eta^{\ell} (x) \;=\; \frac{1}{\ell^d}  \sum_{y\in \bb B_+^\ell(x)} \eta(y)\;.
\end{equation*}

Fix $T>0$ and let $D([0,T], \mc M)$ be the space of $\mc M$-valued
c\`adl\`ag trajectories $\pi:[0,T]\to\mc M$ endowed with the
\emph{uniform} topology.  For each probability measure $\mu_N$ on
$\{0,1\}^{\bb T^d_N}$, denote by $\bb Q_{\mu_N}^{W,N}$ the measure on
the path space $D([0,T], \mc M)$ induced by the measure $\mu_N$ and
the process $\pi^N_t$ introduced in \eqref{f01}.

Fix a continuous profile $\rho_0 : \bb T^d \to [0,1]$ and consider a
sequence $\{\mu_N : N\ge 1\}$ of measures on $\{0,1\}^{\bb T^d_N}$
associated to $\rho_0$ in the sense \eqref{f09}. Further, we denote by $\bb Q_{W}$ be
the probability measure on $D([0,T], \mc M)$ concentrated on the
deterministic path $\pi(t,du) = \rho (t,u)du$, where $\rho$ is the
unique weak solution of \eqref{g03} with $\gamma = \rho_0$, $l_k=0$,
$r_k=1$, $k=1,\ldots,d$ and $\Phi(\alpha) = \alpha + a\alpha^2$.

In subsection \ref{ss1} we show that the sequence $\{\bb Q_{\mu_N}^{W,N} : N\ge
1\}$ is tight and in subsection \ref{ss2} we characterize the limit
points of this sequence.

\subsection{Tightness}
\label{ss1}

The proof of tightness of sequence $\{\bb Q_{\mu_N}^{W,N} : N\ge 1\}$ is
motivated by \cite{jl,TC}. We consider initially the auxiliary $\mc
M$-valued Markov process $\{\Pi^{\lambda,N}_t : t\ge 0\}$, $\lambda>0$,
defined by
\begin{equation*}
\Pi^{\lambda,N}_t (H) \;=\; \< \pi^N_t, G_\lambda^N H\>\;=\;
\frac{1}{N^d} \sum _{x\in \bb Z^d} \bigl (G_\lambda^N H \bigr)(x/N)
\eta_t (x),
\end{equation*}
for $H$ in $C(\bb T^d)$, where $\{G_\lambda^N : \lambda >0\}$ is the
resolvent associated to the random walk $\{X^N_t : t\ge 0\}$
introduced in Section \ref{sec04}.

We first prove tightness of the process $\{\Pi^{\lambda,N}_t : 0\le t
\le T\}$ for every $\lambda>0$ and we then show that $\{\Pi^{\lambda,N}_t
: 0\le t \le T\}$, and that $\{\pi^{N}_t : 0\le t \le T\}$ are not far
apart if $\lambda$ is large.

It is well known \cite{kl} that to prove
tightness of $\{\Pi^{\lambda,N}_t : 0\le t \le T\}$ it is enough to
show tightness of the real-valued processes $\{\Pi^{\lambda,N}_t (H) :
0\le t \le T\}$ for a set of smooth functions $H:\bb T^d\to \bb R$ dense
in $C(\bb T^d)$ for the uniform topology.

Fix a smooth function $H: \bb T^d \to \bb R$. Denote by the same symbol
the restriction of $H$ to $N^{-1} \bb T^d_N$. Let $H_\lambda^N =
G_\lambda^N H$, so that
\begin{equation}
\label{f04}
\lambda H_\lambda^N \;-\; \bb L_N H_\lambda^N \;=\; H\;.
\end{equation}
Keep in mind that $\Pi^{\lambda,N}_t (H) = \<\pi^N_t, H_\lambda^N \>$,
and denote by $M^{N,\lambda}_t$ the martingale defined by
\begin{equation}
\label{f10}
M^{N,\lambda}_t \;=\;  \Pi^{\lambda,N}_t (H) \;-\;
\Pi^{\lambda,N}_0 (H) \;-\; \int_0^t ds \, N^2 L_N \<\pi^N_s ,
H_\lambda^N \> \;.
\end{equation}
Clearly, tightness of $\Pi^{\lambda,N}_t (H)$ follows from tightness
of the martingale $M^{N,\lambda}_t$ and tightness of the additive
functional $\int_0^t ds \, N^2 L_N \<\pi^N_s , H_\lambda^N \>$.

A long, but simple, computation shows that the quadratic variation
$\<M^{N,\lambda}\>_t$ of the martingale $M^{N,\lambda}_t$ is given by:
$$ \frac{1}{N^{2d}}\sum^d_{j=1}\sum_{x\in \bb T^d}\xi_{x,x+e_j}[\nabla_{N,j}H^N_{\lambda}(x/N)]^2
\int_0^t c_{x,x+e_j}(\eta_s) \, [\eta_s(x+e_j) - \eta_s(x)]^2 \, ds\;.$$
 In particular, by \eqref{f05},
\begin{equation*}
  \<M^{N,\lambda}\>_t \;\le\; \frac{C_0 t}{N^{2d}} \sum^d_{j=1} \sum_{x\in \bb T^d_N}
  \xi_{x,x+e_j} \, [(\nabla_{N,j} H_\lambda^N)(x/N)]^2 \;\le\; \frac{C(H)t}{\lambda N^d},
\end{equation*}
for some finite constant $C(H)$ which depends only on $H$. Thus, by
Doob inequality, for every $\lambda>0$, $\delta>0$,
\begin{equation}
\label{f02}
\lim_{N\to\infty} \bb P_{\mu_N} \left[ \sup_{0\le t\le T}
\big\vert M^{N,\lambda}_t \big\vert \, > \, \delta \right]
\;=\; 0\;.
\end{equation}
In particular, the sequence of martingales $\{M^{N,\lambda}_t : N\ge
1\}$ is tight for the uniform topology.

It remains to examine the additive functional of the decomposition
\eqref{f10}. The  generator of the exclusion process $L_N$  is decomposed
 in generators of the random walks $\bb L_{N,j}$.  By
 \eqref{ff}, \eqref{fg} and a long but simple computation, we obtain that
 $N^2 L_N \<\pi^N , H_\lambda^N \>$ is equal to
\begin{eqnarray*}
\!\!\!\!\!\!\!\!\!\!\!\!\!\! &&
\sum^d_{j=1}\big \{\frac {1}{N^d} \sum_{x\in \bb T^d_N} (\bb L_N^j H_\lambda^N)(x/N)\, \eta(x)
\\
\!\!\!\!\!\!\!\!\!\!\!\!\!\! && \quad
+\; \frac{a}{N^d} \sum_{x\in \bb T^d_N} \big [ (\bb L_N^j H_\lambda^N)
((x+e_j)/N) + (\bb L_{N,j} H_\lambda^N) (x/N) \big ] \,
(\tau_x h_{1,j}) (\eta) \\
\!\!\!\!\!\!\!\!\!\!\!\!\!\! && \qquad
- \; \frac{a}{N^d} \sum_{x\in \bb T^d_N} (\bb L_N^j H_\lambda^N)
(x/N) (\tau_x h_{2,j}) (\eta)\big \}\;,
\end{eqnarray*}
where $\{\tau_x: x\in \bb Z^d\}$ is the group of translations, so that
$(\tau_x \eta)(y) = \eta(x+y)$ for $x$, $y$ in $\bb Z^d$, and the
sum is understood modulo $N$. Also, $h_{1,j}$, $h_{2,j}$ are the cylinder functions
\begin{equation*}
h_{1,j}(\eta) \;=\; \eta(0) \eta({e_j})\;,\quad h_{2,j}(\eta) \;=\; \eta(-e_j)  \eta(e_j)\;.
\end{equation*}

Since $H_\lambda^N$ is the solution of the resolvent equation
\eqref{f04}, we may replace $\bb L_N H_\lambda^N$ by $U_\lambda^N
= \lambda H_\lambda^N - H$ in the previous formula. In particular,
for all $0\le s<t\le T$,
\begin{equation*}
\Big\vert \int_s^t dr \, N^2 L_N \<\pi^N_r ,H_\lambda^N \> \Big\vert
\;\le\; \frac {(1+3|a|)(t-s)}{N^d} \sum_{x\in \bb T^d_N} |U_\lambda^N (x/N)|
\;.
\end{equation*}
It follows from the first estimate in \eqref{f05}, and from Schwarz
inequality, that the right hand side of the previous expression  is bounded above by $C(H,a) (t-s)$
uniformly in $N$, where $C(H,a)$ is a finite constant depending only
on $a$ and $H$. This proves that the additive part of the
decomposition \eqref{f10} is tight for the uniform topology and
therefore that the sequence of processes $\{\Pi^{\lambda,N}_t :N\ge
1\}$ is tight.

\begin{lemma}
\label{s06}
The sequence of measures $\{\bb Q_{\mu^N}^{W,N} : N\ge 1\}$ is tight
for the uniform topology.
\end{lemma}

\begin{proof}
It is enough to show that for every smooth function $H:\bb T\to\bb R$
and every $\epsilon>0$, there exists $\lambda>0$ such that
\begin{equation*}
\lim_{N\to\infty} \bb P_{\mu^N} \left[
\sup_{0\le t\le T} |\, \Pi^{\lambda,N}_t (\lambda H) -
\<\pi^N_t, H\>\, | > \epsilon
\right] \;=\;0,
\end{equation*}
since in this case,  the tightness of $\pi^N_t$ follows from 
tightness of $\Pi^{\lambda,N}_t$.  Since there is at most one particle
per site, the expression inside the absolute value is less than or
equal to
\begin{equation*}
\frac{1}{N^d} \sum_{x \in \bb T^d_N} \big|\lambda H_\lambda^N (x/N)
- H (x/N)\big|\;.
\end{equation*}
By Lemma \ref{2.8} this expression vanishes as $N\uparrow\infty$ and then
$\lambda\uparrow\infty$.
\end{proof}

\subsection{Uniqueness of limit points}
\label{ss2}

We prove in this subsection that all limit points $\bb Q^*$ of the
sequence $\bb Q^{W,N}_{\mu_N}$ are concentrated on absolutely
continuous trajectories $\pi(t,du) = \rho(t,u) du$, whose density
$\rho(t,u)$ is a weak solution of the hydrodynamic equation
\eqref{g03} with $l=0$ $<r=1$ and $\Phi(\alpha)=\alpha + a\alpha^2$.

Let $\bb Q^*$ be a limit point of the sequence $\bb Q^{W,N}_{\mu_N}$
and assume, without loss of generality, that $\bb Q^{W,N}_{\mu_N}$
converges to $\bb Q^*$.

Since there is at most one particle per site, it is clear that $\bb
Q^*$ is concentrated on trajectories $\pi_t(du)$ which are absolutely
continuous with respect to the Lebesgue measure, $\pi_t(du) =
\rho(t,u) du$, and whose density $\rho$ is non-negative and bounded by
$1$.

Fix a continuously differentiable function $H: \bb T^d\to \bb R$  and
$\lambda>0$.  Recall the definition of the martingale
$M^{N,\lambda}_t$ introduced in the previous section. By \eqref{f02},
for every $\delta>0$,
\begin{equation*}
\lim_{N\to\infty} \bb P_{\mu_N} \left[ \sup_{0\le t\le T}
\big\vert M^{N,\lambda}_t \big\vert \, > \, \delta \right]
\;=\; 0\;.
\end{equation*}

By \eqref{f10}, for fixed $0<t\le T$ and $\delta>0$,
\begin{equation*}
\lim_{N\to\infty} \bb Q^{W,N}_{\mu_N} \left[ \,
\Big\vert \<\pi^N_t, G^N_\lambda H \> \;-\;
\<\pi^N_0, G^N_\lambda H \> \;-\;
\int_0^t ds \, N^2 L_N \<\pi^N_s , G^N_\lambda H \>
\Big\vert \, > \, \delta \right] \;=\; 0.
\end{equation*}

Since there is at most one particle per site, we may
replace, by Corollary \ref{c1}, $G^N_\lambda H$ by $G_\lambda H$ in the expressions $\<\pi^N_t,
G^N_\lambda H \>$, $\<\pi^N_0, G^N_\lambda H \>$ above.
On the other hand, the expression $N^2 L_N \<\pi^N_s , G^N_\lambda H
\>$ has been computed in the previous subsection. Recall that $\bb L_N
G^N_\lambda H = \lambda G^N_\lambda H - H$. As before, we may replace
$G^N_\lambda H$ by $G_\lambda H$. Let $U_\lambda = \lambda G_\lambda H
- H$.  Since $E_{\nu_\alpha}[h_{i,j}] = \alpha^2$, $i=1$, $2$ and $j = 1,\ldots, d$, in view of
\eqref{f05}, and by Corollary \ref{s02}, for every $t>0$, $\lambda>0$,
$\delta>0$, $i=1$, $2$,
\begin{equation*}
\lim_{\varepsilon \to 0} \limsup_{N\to\infty}
\bb P_{\mu_N} \left[ \, \Big| \int_0^t \!\!\! ds\, \frac 1{N^d}
\sum_{x\in \bb T^d_N} U_\lambda  (x/N) \left\{ \tau_x h_{i,j} (\eta_s) -
\left[\eta^{\varepsilon N}_s(x)\right]^2 \right\} \, \Big|
\, > \, \delta \, \right]  \;=\; 0.
\end{equation*}

Since $\eta^{\varepsilon N}_s(x) = \varepsilon^{-d} \pi^N_s (\prod_{j=1}^d[x_j/N, x_j/N + \varepsilon e_j])$,
 we obtain, from the previous considerations, that
\begin{align*}
\lim_{\varepsilon \to 0} \limsup_{N\to\infty} \bb Q^{W,N}_{\mu_N} \left[ \,
\Big\vert\right. &\<\pi^N_t, G_\lambda H \> \;-\; \\
-\; \<\pi^N_0, G_\lambda H \> \;-\;&\left.
\int_0^t ds \, \Big\< \Phi \big(\varepsilon^{-d} \pi^N_s (\prod_{j=1}^d[\cdot, \cdot
+ \varepsilon e_j]) \big) \,,\, U_\lambda\Big>
\Big\vert > \delta \right] \;=\; 0\;.
\end{align*}

Since $H$ is a smooth function, $G_\lambda H$ and $U_\lambda$ can be
approximated, in $L^1(\bb T^d)$, by continuous functions. Since we assumed
that $\bb Q^{W,N}_{\mu_N}$ converges in the uniform topology to $\bb
Q^*$, we have that
\begin{align*}
\lim_{\varepsilon \to 0}\bb Q*\left[ \,
\Big\vert \<\pi_t, G_\lambda H \> \right.&\;-\; \; \<\pi_0, G_\lambda H \> \;-\;\\
-\;\int_0^t ds \, &\left.\Big\< \Phi \big (\varepsilon^{-d} \pi_s (\prod_{j=1}^d[\cdot, \cdot
+ \varepsilon e_j]) \big) \,,\, U_\lambda\Big>
\Big\vert > \delta \right] \;=\; 0\;.
\end{align*}

Since $\bb Q^{*}$ is concentrated on absolutely continuous paths
$\pi_t(du) = \rho(t,u) du$ with positive density bounded by $1$,
$\varepsilon^{-d}\pi_s(\prod_{j=1}^d[\cdot, \cdot + \varepsilon e_j])$ converges in
$L^1(\bb T^d)$ to $\rho(s,.)$ as $\varepsilon\downarrow 0$. Thus,
\begin{eqnarray*}
\bb Q^{*} \left[ \,
\Big\vert \<\pi_t, G_\lambda H \> \;-\;
 \<\pi_0, G_\lambda H \> \;-\;
\int_0^t ds \, \< \Phi (\rho_s) \,,\, \mc L_W G_\lambda H \>
\Big\vert > \delta \right] \;=\; 0,
\end{eqnarray*}
because $U_\lambda = \mc L_W G_\lambda H$. Letting $\delta\downarrow
0$, we see that, $\bb Q^{*}$ a.s.,
\begin{eqnarray*}
\<\pi_t, G_\lambda H \> \;-\; \<\pi_0, G_\lambda H \> \;=\;
\int_0^t ds \, \< \Phi (\rho_s) \,,\, \mc L_W G_\lambda H \> \;.
\end{eqnarray*}
This identity can be extended to a countable set of times $t$. Taking
this set to be dense, by continuity of the trajectories $\pi_t$, we
obtain that it holds for all $0\le t\le T$. In the same way, it holds
for any countable family of continuous functions $H$. Taking a countable
set of continuous functions, dense for the uniform topology, we extend
this identity to all continuous functions $H$, because $G_\lambda H_n$
converges to $G_\lambda H$ in $L^1(\bb T^d)$, if $H_n$ converges to $H$
in the uniform topology.  Similarly, we can show that it holds for all
$\lambda>0$, since, for any continuous function $H$, $G_{\lambda_n} H$
converges to $G_\lambda H$ in $L^1(\bb T^d)$, as $\lambda_n \to
\lambda$.

\begin{proposition}
\label{s15}
As $N\uparrow\infty$, the sequence of probability measures $\bb
Q_{\mu_N}^{W,N}$ converges in the uniform topology to $\bb Q_{W}$.
\end{proposition}
\begin{proof}
In the previous subsection we showed that the sequence of probability
measures $\bb Q^{W,N}_{\mu_N}$ is tight for the uniform topology. Moreover, we
just proved that all limit points of this sequence are concentrated on
weak solutions of the parabolic equation \eqref{g03}. The proposition now follows
 from a straightforward adaptation of the   uniqueness of weak solutions 
 proved in \cite{TC} for the $d$-dimensional case.
\end{proof}

\begin{proof}[Proof of Theorem \ref{t02}]
Since $\bb Q_{\mu_N}^{W,N}$ converges in the uniform topology to $\bb
Q_{W}$, a measure which is concentrated on a deterministic path, for
each $0\le t\le T$ and each continuous function $H:\bb T^d\to \bb R$,
$\<\pi^N_t, H\>$ converges in probability to $\int_{\bb T} du \,
\rho(t,u)$ $H(u)$, where $\rho$ is the unique weak solution of
\eqref{g03} with $l_k=0$, $r_k=1$, $\gamma=\rho_0$ and $\Phi(\alpha) =
\alpha + a \alpha^2$.
\end{proof}

\subsection{Replacement lemma}

We will use some results from \cite[Appendix A1]{kl}.
 Denote by $H_N (\mu_N | \nu_\alpha)$ the relative entropy of a probability
measure $\mu_N$ with respect to a stationary state $\nu_\alpha$,
see \cite[Section A1.8]{kl} for a precise definition. By the
explicit formula given in \cite[Theorem A1.8.3]{kl}, we see that there
exists a finite constant $K_0$, depending only on $\alpha$, such that
\begin{equation}
\label{f06}
H_N (\mu_N | \nu_\alpha) \;\le\; K_0 N^d,
\end{equation}
for all measures $\mu_N$.

Denote by $\< \cdot, \cdot \>_{\nu_\alpha}$ the scalar product of
$L^2(\nu_\alpha)$ and denote by $I^\xi_N$ the convex and lower
semicontinuous \cite[Corollary A1.10.3]{kl} functional defined by
\begin{equation*}
I^\xi_N (f) \;=\; \< - L_N \sqrt f \,,\, \sqrt f\>_{\nu_\alpha}\; ,
\end{equation*}
for all probability densities $f$ with respect to $\nu_\alpha$ (i.e.,
$f\ge 0$ and $\int f d\nu_\alpha =1$). By \cite[ proposition A1.10.1]{kl},
an elementary computation shows  that
\begin{eqnarray*}
\!\!\!\!\!\!\!\!\!\!\!\!\!\! &&
I^\xi_N (f) \;=\;\sum^d_{j=1} \sum_{x\in \bb T^d_N} I^\xi_{x,x+e_j} (f)\;,
\quad\text{where}\quad \\
\!\!\!\!\!\!\!\!\!\!\!\!\!\! && \qquad
I^\xi_{x,x+e_j} (f) \;=\; (1/2)\, \xi_{x,x+e_j}
\int c_{x,x+e_j} (\eta) \left\{ \sqrt{f(\sigma^{x,x+e_j} \eta)} -
\sqrt{f(\eta)} \right\}^2 \, d\nu_\alpha \;.
\end{eqnarray*}
By \cite[Theorem A1.9.2]{kl}, if $\{S^N_t : t\ge 0\}$ stands for the
semi-group associated to the generator $N^2L_N$,
\begin{equation*}
H_N (\mu_N S^N_t | \nu_\alpha) \; + 2\; N^2 \, \int_0^t
I^\xi_N (f^N_s) \, ds  \;\le\; H_N (\mu_N | \nu_\alpha)\;,
\end{equation*}
where $f^N_s$ stands for the Radon-Nikodym derivative of $\mu_N
S^N_s$ with respect to $\nu_\alpha$.
\medskip

Remember that for a local function $g: \{0,1\}^{\bb Z^d} \to \bb R$, 
$\tilde g :[0,1]\to \bb R$ stands for the expected value of $g$ under the
stationary states:
\begin{equation*}
\tilde g (\alpha) \;=\; E_{\nu_\alpha} [ g(\eta)]\;.
\end{equation*}
For $\ell \ge 1$ and $d$-dimensional integer $x=(x_1,\ldots,x_d)$, denote 
by $\eta^{\ell} (x)$ the empirical density of particles in the box 
$ \bb B_+^\ell(x)= \{(y_1,\ldots,y_d)\in\bb Z^d\;;0\le y_i-x_i < \ell\}$: 
\begin{equation*}
  \eta^{\ell} (x) \;=\; \frac{1}{\ell^d}  \sum_{y\in \bb B_+^\ell(x)} \eta(y)\;.
\end{equation*}

For each $y\in \bb B^\ell_+(x)$, such that $y_1>x_1$, let 
\begin{equation}
\label{ca}
\Lambda^\ell_{x+e_1,y}= (z^y_k)_{0\le k\le M(y)}
\end{equation}
be a path from $x+e_1$ to $y$ such that:
\begin{enumerate}
\item $\Lambda^\ell_{x+e_1,y}$ begins at $x+e_1$ and ends at $y$, i.e.:
$$z^y_0\;=\; x+e_1 \;\;\;\text{and}\;\;\;z^y_{M(y)}=y;$$
\item The distance between two consecutive sites of the $\Lambda^\ell_{x+e_1,y}= (z^y_k)_{0\le k\le M(y)}$ 
is equal  to $1$, i.e.:
\begin{equation*}
z^y_{k+1}=z^y_k +e_j ;\;\text{for some  }j=1\ldots,d \;\;\;\text{and for all} \;\;k=1,\ldots, M(y)-1
\end{equation*}
\item The number of points $M(y)$ is bounded above by $d\ell$;
\item $\Lambda^\ell_{x+e_1,y}$ is injective:
$$z^y_i\neq z^y_j\;\;\; \text{for all  } 0\le i < j \le M(y).$$
\end{enumerate}
\begin{lemma}
\label{s01}
Fix a function $F: N^{-1} \bb T^d_N \to \bb R$. There exists a finite constant
$C_0 = C_0(a,g,W)$, depending only on $a$, $g$ and $W$, such that
\begin{eqnarray*}
\!\!\!\!\!\!\!\!\!\!\!\! &&
\frac{1}{N^d} \sum_{x\in \bb T^d_N} F(x/N) \int  \{ \tau_x g (\eta) - \tilde
g (\eta^{\varepsilon N}(x)) \}\,  f(\eta) \nu_\alpha(d\eta) \\
\!\!\!\!\!\!\!\!\!\!\!\! && \quad
\le\; \frac{C_0}{\varepsilon N^{d+1}} \sum_{x\in \bb T_N^d} \big| F(x/N) \big|
\;+\; \frac{C_0 \varepsilon}{\delta N^d} \sum_{x\in \bb T^d_N} F(x/N)^2
 \;+\; \frac {\delta}{N^{d-2}} I^\xi_N(f),
\end{eqnarray*}
for all $\delta>0$, $\varepsilon >0$ and all probability densities $f$ with respect to
$\nu_\alpha$.
\end{lemma}

\begin{proof}
Any local function can be written as a linear combination of functions
of type $\prod_{x\in A} \eta(x)$, for finite sets $A's$. It is
therefore enough to prove the lemma for such functions. We will only prove the
result for $g(\eta) = \eta(0) \eta(e_1)$. The general case can be
handled in a similar way.

We begin by estimating
\begin{equation}
\label{f03}
\frac {1}{N^d} \sum_{x\in \bb T^d_N} F(x/N) \int \eta(x) \{ \eta(x+e_1) -
\frac 1{(\varepsilon N)^d}  \sum_{y\in \bb B^{N\varepsilon}_+(x)} \eta(y)\} f(\eta) \nu_\alpha(d\eta)
\end{equation}
in terms of the functional $I^\xi_N(f)$. 
The integral in \eqref{f03} can be rewritten as:
\begin{equation*}
\frac{1}{(N\varepsilon)^d}\sum_{y\in \bb B^{N\varepsilon}_+(x)}\int \eta(x)[\eta(x+e_1)-\eta(y)]f(\eta)\nu_\alpha(d\eta)
\end{equation*}

For each $y\in \bb B^{N\varepsilon}_+(x)$, such that $y_1>x_1$, let $\Lambda^\ell_{x+e_1,y}= (z_k^y)_{0\le k\le M(y)}$
be a path like the one in \eqref{ca}. Then, by property $(1)$ of $\Lambda^\ell_{x+e_1,y}$ and using telescopic sum we have the following:
\begin{equation*}
\eta(x+e_1)-\eta(y) = \sum_{k=0}^{M(y)-1}[\eta(z^y_k)-\eta(z^y_{k+1})].
\end{equation*}

We can, therefore, bound (\ref{f03}) above by 
\begin{align*}
\frac {1}{N^d} \frac{1}{(N\varepsilon)^d} \sum_{x\in \bb T^d_N}
&\sum_{y\in \bb B^{N\varepsilon}_+(x)}\sum_{k=0}^{M(y)-1}
\int F(x/N) \eta(x)[\eta(z^y_k)-\eta(z^y_{k+1})]f(\eta)\nu_\alpha(d\eta)\;+\\
&\frac{1}{\varepsilon N^{d+1}} \sum_{x\in \bb T^d_N} \big| F(x/N) \big|
\end{align*}
where the last term in the previous expression comes from the contribution of the points
$y\in \bb B^{N\varepsilon}_+(x)$, such that $y_1=x_1$.
Recall that by property $(2)$ of $\Lambda^\ell_{x+e_1,y}$, we have that $z^y_{k+1} = z^y_{k} +e_j$,
for some $j=1,\ldots,d$.

For each term of the form 
\begin{equation*}
\int F(x/N)\eta(x)\{\eta(z)-\eta(z+e_j)\}f(\eta)\nu_\alpha(d\eta)
\end{equation*}
we can use the change of
variables $\eta' = \sigma^{z ,z + e_j}\eta$ to write the previous integral as
\begin{equation*}
(1/2) \int F(x/N)\eta(x) \{ \eta(z) - \eta(z+e_j) \} \,
\left\{ f(\eta) - f(\sigma^{z,z+e_j}\eta) \right\} \, \nu_\alpha(d\eta) \;.
\end{equation*}
Since $a-b = (\sqrt a - \sqrt b)(\sqrt a + \sqrt b)$ and $\sqrt{ab}\le a+b$, by Schwarz
inequality the previous expression is less than or equal to
\begin{align*}
\frac A {4 (1-2a^-) \xi_{z,z+e_j}} \int F(x/N)^2\eta(x)& \{ \eta(z) - \eta(z+e_j)\}^2 \times\\
\times& \left\{ \sqrt{f(\eta)} + \sqrt{f(\sigma^{z,z+e_j}\eta)} \right\}^2 \, \nu_\alpha(d\eta) \;+\;\\
+\frac {\xi_{z,z+e_j}} A  \int c_{z,z+e_j}(\eta)&
\left\{ \sqrt{f(\eta)} - \sqrt{f(\sigma^{z,z+e_j}\eta)} \right\}^2
\, \nu_\alpha(d\eta)
\end{align*}
for every $A>0$. In this formula we used the fact that $c_{z,z+e_j}$ is
bounded below by $1-2a^-$. Since $f$ is a density with respect to
$\nu_\alpha$, the first expression is bounded above by $A/(1-2a^-) \xi_{z,z+e_j}$,
whereas the second one is equal to $2 A^{-1} I^\xi_{z,z+e_j}(f)$. 

So, by properties $(3)$ and $(4)$ of the path $\Lambda^\ell_{x+e_1,y}$, we obtain that \eqref{f03} is less than or equal to
\begin{eqnarray*}
\!\!\!\!\!\!\!\!\!\!\!\!\!\! &&
\frac{1}{\varepsilon N^{d+1}} \sum_{x\in \bb T^d_N} \big| F(x/N) \big|
\;+\; \frac A{(1-2a^-) N^d} \sum_{x\in \bb T^d_N} F(x/N)^2
\sum_{j=1}^d\sum_{k=1}^{\varepsilon N} \xi_{x+(k-1)e_j,x+ke_j}^{-1}
\;+\; \\
\!\!\!\!\!\!\!\!\!\!\!\!\!\! &&
\;\;\;\;\;\;\;\;\;\frac {2\varepsilon} {AN^{d-1}}\sum_{j=1}^d\sum_{x\in \bb T^d_N} I^\xi_{x, x+e_j}(f)\;.
\end{eqnarray*}
By definition of the sequence $\{\xi_{x,x+e_j}\}$, $\sum_{k=1}^{\varepsilon N}
 \xi_{x+ke_j,e_j}^{-1} \le N [W_j(1) - W_j(0)]$.
Thus, choosing $A=2 \varepsilon N^{-1} \delta^{-1}$, for some $\delta>0$,
we obtain that the previous sum is bounded above by
\begin{equation*}
\frac {C_0}{\varepsilon N^{d+1}} \sum_{x\in \bb T^d_N} \big| F(x/N) \big|
\;+\; \frac {C_0 \varepsilon}{\delta N^d} \sum_{x\in \bb T^d_N} F(x/N)^2
 \;+\; \frac {\delta}{N^{d-2}} I^\xi_N(f)\;.
\end{equation*}

Up to this point we have succeeded to replace $\eta(x) \eta(x+e_1)$ by $\eta(x)
\eta^{\varepsilon N} (x)$. The same arguments permit to replace this
latter expression by $[\eta^{\varepsilon N} (x)]^2$, which concludes
the proof of the lemma.
\end{proof}

\begin{corollary}
\label{s02}
Fix a cylinder function $g$ and a sequence of functions $\{F_N : N\ge
1\}$, $F_N : N^{-1} \bb T^d_N\to \bb R$ such that
\begin{equation*}
\limsup_{N\to\infty} \frac 1{N^d} \sum_{x\in \bb T^d_N} F_N(x/N)^2 \;<\;
\infty\; .
\end{equation*}
Then, for any $t>0$ and any sequence of probability measures $\{\mu_N
: N\ge 1\}$ on $\{0,1\}^{\bb T^d_N}$,
\begin{equation*}
\limsup_{\varepsilon\to 0} \limsup_{N\to\infty}
\bb E_{\mu_N} \Big[ \, \Big| \int_0^t  \frac 1{N^d}
\sum_{x\in \bb T_N^d} F_N(x/N) \, \big \{ \tau_x g (\eta_s) - \tilde g
(\eta^{\varepsilon N}_s(x))\ d_s \big \} \Big| \, \Big]  \;=\; 0\;.
\end{equation*}
\end{corollary}

\begin{proof}
Fix $0<\alpha <1$. By the entropy and Jensen inequalities, the
expectation appearing in the statement of the lemma is bounded above
by
$$\frac 1{\gamma N^d} \log \bb E_{\nu_\alpha} \Big[
\exp\Big\{ \gamma \Big| \int_0^t \!\!\! ds \!\!\!
\sum_{x\in \bb T^d_N} F_N(x/N) \, \big \{ \tau_x g (\eta_s) - \tilde g
(\eta^{\varepsilon N}_s(x)) \big \} \, \Big| \, \Big\} \, \Big]\\
+\frac {H_N (\mu_N | \nu_\alpha)}{\gamma N^d} 
$$
for all $\gamma >0$. In view of \eqref{f06}, in order to prove the corollary it
is enough to show that the second term vanishes as $N\uparrow\infty$,
and then $\varepsilon\downarrow 0$ for every $\gamma>0$. We may remove
the absolute value inside the exponential by using the elementary inequalities $e^{|x|} \le e^x +
e^{-x}$ and $\limsup_{N\to\infty} N^{-1} \log\{a_N + b_N\} \le
\max\{ \limsup_{N\to\infty} N^{-1} \log a_N , \limsup_{N\to\infty}
N^{-1}$ $\log b_N \}$. Thus, to prove the corollary, it is enough to show
that
\begin{equation*}
\limsup_{\varepsilon\to 0} \limsup_{N\to\infty}
\frac 1{N^d} \log \bb E_{\nu_\alpha} \Big[
\exp\Big\{ \gamma \int_0^t \!\!\! ds \!\!\!
\sum_{x\in \bb T^d_N} F_N(x/N) \{ \tau_x g (\eta_s) - \tilde g
(\eta^{\varepsilon N}_s(x)) \} \Big\} \, \Big] = 0
\end{equation*}
for every $\gamma>0$.

By Feynman-Kac formula, for each fixed $N$ the previous expression is
bounded above by
\begin{equation*}
t\gamma\, \sup_{f} \left\{ \int  \frac 1 {N^d}
\sum_{x\in \bb T^d_N} F_N(x/N) \{ \tau_x g (\eta) - \tilde g
(\eta^{\varepsilon N}(x)) \} f (\eta) \, d \nu_\alpha\;
-\; \frac{1}{N^{d-2}}  \hat I^\xi_N (f) \right\},
\end{equation*}
where  the supremum  is carried  over all  density functions  $f$ with
respect to  $\nu_\alpha$. Letting $\delta  =1$ in Lemma  \ref{s01}, we
obtain that the previous expression is less than or equal to
\begin{equation*}
\frac {C_0\gamma t} {\varepsilon N^{d+1}}
\sum_{x\in \bb T^d_N} \big| F_N(x/N) \big|
\;+\; \frac {C_0\gamma \varepsilon t}{N^d}
\sum_{x\in \bb T^d_N} F_N(x/N)^2,
\end{equation*}
for some finite constant $C_0$ which depends on $a$, $g$ and $W$. By
assumption on the sequence $\{F_N\}$, for every $\gamma>0$, this
expression vanishes as $N\uparrow\infty$ and then
$\varepsilon\downarrow 0$. This concludes the proof of the lemma.
\end{proof}

\section{Energy estimate}
\label{sec4}

We prove in this section that any limit point $\bb Q^*_{W}$ of the
sequence $\bb Q_{\mu_N}^{W,N}$ is concentrated on trajectories
$\rho(t,u) du$ with finite energy.

Denote by
$\partial_{x_j}$ the partial derivative of a function with respect to the
$j$-th coordinate, and by $C^{0,1_j}([0,T]\times \bb T^d)$
the set of continuous functions with continuous partial derivative in the $j$-th coordinate.
Let $L^2_{x_j\otimes W_j}([0,T]\times \bb T^d)$ be the Hilbert space of
measurable functions $H: [0,T]\times \bb T^d\to\bb R$ such that
\begin{equation*}
\int_0^Tds\int_{\bb T^d}d(x_j\otimes W_j) \, H (s, u)^2 \;<\; \infty\;,
\end{equation*}
where $d(x_j\otimes W_j)$ represents the product measure in $\bb T^d$
obtained from Lesbegue's measure in $\bb T^{d-1}$ and the measure induced by $W_j$:
$$d(x_j\otimes W_j)\;=\;dx_1\ldots dx_{j-1}\;dW_j\; dx_{j+1}\ldots dx_d\;,$$
endowed with the inner product $\<\!\< H,G \>\!\>_{x_j\otimes W_j}$ defined by
\begin{equation*}
\<\!\< H,G \>\!\>_{x_j\otimes W_j} \;=\; \int_0^T ds \int_{\bb T^d} d(x_j\otimes W_j) \, H (s, u)
\, G(s,u)\;.
\end{equation*}

Let $\bb Q^*_{W}$ be a limit point of the sequence $\bb
Q_{\mu_N}^{W,N}$ and assume without loss of generality that the
sequence $\bb Q_{\mu_N}^{W,N}$ converges to $\bb Q^*_{W}$.  
\begin{proposition}
\label{s05}
The measure $\bb Q^*_{W}$ is concentrated on paths $\rho(t,x) dx$ with
the property that for all $j=1,\ldots,d$ there exists a function in $L^2_{x_j\otimes W_j}([0,T]\times \bb T^d)$,
denoted by $d\Phi/dW_j$, such that
\begin{align*}
\int_0^Tds \int_{\bb T^d} dx& \, (\partial_{x_j} H) (s, x) \, \Phi(\rho(s,x))
\;=\;\\ -\; &\int_0^Tds \int_{\bb T} d(x_j\otimes W_j(x)) \, (d\Phi/dW_j) (s, x) \, H (s, x)
\end{align*}
for all functions $H$ in $C^{0,1_j}([0,T]\times \bb T^d)$.
\end{proposition}

The previous proposition follows from the next lemma.  Recall the
definition of the constant $K_0$ given in \eqref{f06}.

\begin{lemma}
\label{s03}
There exists a finite constant $K_1$, depending only on $a$, such that
\begin{align*}
E_{\bb Q^*_{W}} \left[ \sup_H \left\{ \int_0^T ds\, \int_{\bb T^d}
dx\right.\right. &\, (\partial_{x_j} H) (s, x) \, \Phi(\rho(s,x)) \\ 
- \; K_1 \int_0^T ds\, &\left.\left. \int_{\bb T^d} H (s, x)^2
\, d(x_j\otimes W_j(x)) \right\} \right] \; \le \; K_0 \; ,
\end{align*}
where the supremum is carried over all functions $H \in
C^{0,1_j}([0,T]\times \bb T^d)$.
\end{lemma}

\begin{proof}[Proof of Proposition \ref{s05}]
Denote by $\ell : C^{0,1_j}([0,T]\times \bb T^d) \to \bb R$ the linear
functional defined by
\begin{equation*}
\ell (H) \;=\; \int_0^T ds\, \int_{\bb T^d}
dx \, (\partial_{x_j} H) (s, x) \, \Phi(\rho(s,x))\;.
\end{equation*}
Since $C^{0,1}([0,T]\times \bb T^d)$ is dense in $L^2_{x_j\otimes W_j}([0,T]\times \bb T^d)$,
by Lemma \ref{s03}, $\ell$ is $\bb Q^*_{W}$-almost surely finite
in $L^2_{x_j\otimes W_j}([0,T]\times \bb T^d)$. In particular, by Riesz representation
theorem, there exists a function $G$ in $L^2_{x_j\otimes W_j}([0,T]\times \bb T^d)$
such that
\begin{equation*}
\ell (H) \;=\; - \int_0^T ds\, \int_{\bb T^d}
d(x_j\otimes W_j(x)) \, H (s, x) \, G(s,x)\;.
\end{equation*}
This concludes the proof of the proposition.
\end{proof}

For a smooth function $H\colon \bb T^d\to \bb R$, $\delta >0$, $\varepsilon
>0$ and a positive integer $N$, define $W^j_N(\varepsilon, \delta, H,
\eta)$ by
\begin{eqnarray*}
W^j_N(\varepsilon, \delta, H, \eta ) &=&
\sum_{x\in\bb T^d_N} H(x/N) \frac 1{\varepsilon N}
\, \left\{ \Phi(\eta^{\delta N} (x)) -
\Phi(\eta^{\delta N} (x + \varepsilon Ne_j)) \right\} \\
&-& \frac {K_1}{\varepsilon N} \sum_{x\in\bb T^d_N}
H(x/N)^2 \{ W_j([x_j + \varepsilon N +1]/N) - W_j(x_j/N) \}\; .
\end{eqnarray*}

The proof of Lemma \ref{s03} relies on the following result.
\begin{lemma}
\label{s04}
Consider a sequence $\{H_\ell,\, \ell\ge 1\}$ dense in
$C^{0,1}([0,T]\times \bb T^d)$.  For every $k\ge 1$, and every
$\varepsilon >0$,
\begin{equation*}
\limsup_{\delta\to 0} \limsup_{N\to\infty}
\bb E_{\mu^N} \left[ \max_{1\le i\le k} \left\{
\int_0^T W^j_N(\varepsilon, \delta, H_i (s, \cdot) , \eta_s ) \,
ds \right\} \right] \;\le\; K_0\; .
\end{equation*}
\end{lemma}

\begin{proof}
It follows from the replacement lemma that in order to prove
the Lemma we just need to show that
\begin{equation*}
\limsup_{N\to\infty} \bb E_{\mu^N} \left[ \max_{1\le i\le k} \left\{
\int_0^T W^j_N(\varepsilon,  H_i (s, \cdot) , \eta_s ) \, ds \right\}
\right] \;\le\; K_0\; ,
\end{equation*}
where
\begin{eqnarray*}
W^j_N(\varepsilon, H , \eta ) &=&
\frac 1{\varepsilon N} \sum_{x\in\bb T^d_N} H(x/N)
\left\{ \tau_x g(\eta) -  \tau_{x + \varepsilon Ne_j} g(\eta)\right\} \\
&-&  \frac {K_1}{\varepsilon N} \sum_{x\in\bb T^d_N} H(x/N)^2
\{ W_j([x_j+\varepsilon N +1]/N) - W_j(x_j/N) \}\; ,
\end{eqnarray*}
and $g(\eta) = \eta(0) + a \eta(0)\eta(e_j)$.

By the entropy and Jensen's inequalities,
for each fixed $N$, the previous expectation is bounded above by
\begin{equation*}
\frac {H(\mu^N \vert \nu_{\alpha})}{ N^d} \; +\; \frac 1{N^d}
\log \bb E_{\nu_{\alpha}} \left[ \exp\left\{
\max_{1\le i\le k} \left\{ N^d \int_0^T ds\,
W^j_N(\varepsilon,  H_i (s, \cdot) , \eta_s ) \right\} \right\} \right] \; .
\end{equation*} 
By \eqref{f06}, the first term is bounded by $K_0$.  Since $\exp\{
\max_{1\le j\le k} a_j \}$ is bounded above by $\sum_{1\le j\le k}
\exp\{a_j\}$, and since $\limsup_N N^{-d} \log \{a_N + b_N\}$ is less
than or equal to the maximum of $\limsup_N N^{-d} \log a_N$ and
$\limsup_N N^{-d} \log b_N$, the limit, as $N\uparrow\infty$, of the
second term in the previous expression is less than or equal to
\begin{equation*}
\max_{1\le i \le k} \limsup_{N\to\infty} \frac 1{N^d} \log
\bb E_{\nu_{\alpha}} \left[ \exp
\left\{ N^d  \int_0^T ds\,  W^j_N(\varepsilon,  H_i (s, \cdot) , \eta_s )
\right\} \right] \; .
\end{equation*}
We now prove that, for each fixed $i$, the above limit is non-positive
for a convenient choice of the constant $K_1$.

Fix $1\le i\le k$.
By Feynman--Kac formula and the variational formula for the
largest eigenvalue of a symmetric operator, the previous expression is bounded above by
\begin{equation*}
\int_0^T ds\, \sup_{f} \left\{  \int W^j_N(\varepsilon,
H_i (s, \cdot) , \eta )
f(\eta) \nu_{\alpha} (d\eta) - \frac 1 {N^{d-2}} I^\xi_N (f) \right\}\;,
\end{equation*}
for each fixed $N$.
In this formula the supremum is taken over all probability densities
$f$ with respect to $\nu_{\alpha}$.

To conclude the proof, rewrite $$\eta(x) \eta(x+e_j) - \eta(x+\varepsilon
N e_j)\eta(x+ (\varepsilon N +1)e_j)$$ as $$\eta(x) \{\eta(x+e_j) - \eta(x+
(\varepsilon N +1)e_j)\} + \eta(x+ (\varepsilon N +1)e_j) \{\eta(x) -
\eta(x+ \varepsilon Ne_j)\},$$ and repeat the arguments presented in the
proof of Lemma \ref{s01}.
\end{proof}

\begin{proof}[Proof of Lemma \ref{s03}]
Assume without loss of generality that $\bb Q_{\mu_N}^{W,N}$ converges
to $\bb Q^*_{W}$. Consider a sequence $\{H_\ell,\, \ell\ge 1\}$ dense in
$C^{0,1_j}([0,T]\times \bb T^d)$. By Lemma \ref{s04}, for every $k\ge 1$
\begin{align*}
\limsup_{\delta\to 0} E_{\bb Q^*_{W}}\left[ \max_{1\le i\le k}
\left\{ \frac 1{\varepsilon} \int_0^T ds\, \int_{\bb T^d} dx \,
H_i (s,x) \, \left\{ \Phi(\rho^\delta_s (x)) -
\Phi(\rho^\delta_s (x + \varepsilon e_j)) \right\} \right.\right.\\
-\left.\left.\; \frac {K_1} {\varepsilon} \int_0^T ds\,
\int_{\bb T^d} dx \, H_i(s,x)^2  \, [W_j(x_j+\varepsilon) - W_j(x_j)]
\right\} \right]\; \le \; K_0\; ,
\end{align*}
where $\rho^\delta_s (x) = (\rho_s * \iota_\delta)(x)$ and
$\iota_\delta$ is the approximation of the identity $\iota_\delta
(\cdot) = (\delta)^{-d} \mb 1\{ [0 , \delta]^d\} (\cdot)$.

Letting $\delta\downarrow 0$, changing variables, and then letting
$\varepsilon\downarrow 0$, we obtain that
\begin{align*}
E_{\bb Q^*_{W}}\left[ \max_{1\le i\le k} \left\{
\int_0^T ds\, \int_{\bb T^d} (\partial_{x_j} H_i) (s,x)
\Phi(\rho (s,x)) \, dx \right.\right.\\
-\left.\left. \;K_1 \int_0^T ds\, \int_{\bb T^d} H_i(s,x)^2 d(x_j\otimes W_j(x))
\right\} \right] \;\le \; K_0\; .
\end{align*}
To conclude the proof, it remains to apply the monotone convergence
theorem and recall that $\{H_\ell, \, \ell\ge 1\}$ is a dense sequence
in $C^{0,1_j}([0, T]\times \bb T^d)$ for the norm $\Vert H\Vert_\infty
+ \Vert (\partial_{x_j} H)\Vert_\infty$.
\end{proof}

\section*{Acknowledgments} 
I would like to thank Claudio Landim for giving constant encouragement and several ideas and suggestions that helped in the elaboration of this work. Finally, I would also like to thank Alexandre Simas for his valuable comments.

        \end{document}